\documentclass[11pt]{article}
\usepackage{amsmath}
\usepackage{amssymb}
\usepackage{amsfonts}
\usepackage{theorem}
\usepackage{bbm}

\usepackage[T1]{fontenc}
\usepackage{times}
\usepackage[cp1250]{inputenc}

\newtheorem{dfn}{Definition}[section]
\newtheorem{tw}[dfn]{Theorem}
\newtheorem{prop}[dfn]{Proposition}
\newtheorem{rem}[dfn]{Remark}

\newtheorem{lem}[dfn]{Lemma}

\usepackage{anysize}
\usepackage{verbatim}
\usepackage{array}
\usepackage{enumerate}
\usepackage{amssymb} 
\usepackage{amsmath}
\usepackage{fancyhdr}
\usepackage{euscript}
\usepackage{latexsym}
\usepackage{graphicx}

\author{Micha\l \ Barski \\ \small  Faculty of Mathematics and Computer Science, University of Leipzig, Germany\\
\small Faculty of Mathematics, Cardinal Stefan Wyszy\'nski University in Warsaw, Poland
\\ \small{\it Michal.Barski@math.uni-leipzig.de}}

\title{Approximations for solutions of L\'evy-type stochastic
differential equations \thanks{Research supported by Polish KBN
Grant P03A 034 29 ,,Stochastic evolution equations \ driven by
L\'evy noise''}}
\date{}

\begin{document}
\maketitle \abstract{The problem of the construction of strong
approximations with a given order of convergence for jump-diffusion
equations is studied. General approximation schemes are constructed
for L\'evy type stochastic differential equation. In particular, the
paper generalizes the results from \cite{KloPLa} and \cite{Gar}. The
Euler and the Milstein schemes are shown for finite and infinite
L\'evy measure.}
\\
\begin{quote}

\end{quote}

\section{Introduction}
The problem of approximation construction for solution of stochastic
differential equation is widely studied throughout many papers. The
authors' attention is focused mainly on the equation of the form:
\begin{gather}\label{rownanie ogolne}
Y_t=Y_0+\int_{0}^{t}f(Y_{s-})dZs,
\end{gather}
where $Y_0$ is a random variable with known distribution, $f$-some
regular function and $Z$-a driving process. There are many
approximation methods for the solution of (\ref{rownanie ogolne})
depending on the driving process and the optimality criteria imposed
on the approximating error. The case when $Z$ is a Wiener process
the problem is comprehensively studied in the book \cite{KloPLa},
for jump diffusion case see, for instance, \cite{HigKlo1},
\cite{HigKlo2}. In \cite{KloPLa} various schemes for the so called
weak and strong approximations are presented, in particular their
dependence on the mesh of the partition of the interval $[0,T]$.
Denoting by $\bar{Y}$ the approximation, the optimality criteria for
weak solutions have a form:
$\mathbf{E}[g(Y_T)-g(\bar{Y}_T)]\longrightarrow \min$, where $g$ is
some regular function, while for strong solutions:
$\mathbf{E}\sup_{t}|Y_t-\bar{Y}_t|^2\longrightarrow \min$. The
schemes use the increments of time, increments of the Wiener process
and, for higher order of convergence, some normally distributed
random variables correlated with the increments of the Wiener
process. Thus for practical implementation we have to generate
normally distributed, correlated random variables.

The simplest approximating scheme for the equation $(\ref{rownanie
ogolne})$ is the Euler scheme which has the following structure:
\begin{gather*}
\bar{Y}_0=Y_0,
\qquad\bar{Y}_{\frac{(i+1)T}{n}}=f(\bar{Y}_{\frac{iT}{n}})(Z_{\frac{(i+1)T}{n}}-Z_{\frac{iT}{n}}),
\end{gather*}
where $\{\frac{iT}{n}, i=0,1,...,n\}$ is a partition of the interval
$[0,T]$. In the case of the Wiener driving process it is easy to
construct. However, for a general L\'evy driving process it is no
longer so simple. This is because of the difficulty of practical
construction of the increments of $Z$ when the L\'evy measure is
infinite, i.e. when the measure of a unit ball is infinite. If the
increments can not be simulated, then they themselves have to be
approximated in some sense and the accuracy of such construction
should be studied. This way of approximating is presented for
example in \cite{ProTal} and \cite{Rub}. The main idea in these
papers is to reduce the problem by replacing increments of $Z$ by
suitable increments of the compound Poisson process, which can be
practically simulated. It should be pointed out that our approach is
more general since a significant majority of papers consider
approximation problem using different modifications of the Euler
scheme.

In this paper we work with a stochastic differential equation of the
form:
\begin{align*}
Y_t=Y_0+\int_{0}^{t}b(Y_{s_{-}})ds+\int_{0}^{t}\sigma(Y_{s_{-}})dW_{s}&+\int_{0}^{t}\int_{\mid
x\mid<1}F(Y_{s_{-}},x)\tilde{N}(ds,dx)\\[2ex]
&+\int_{0}^{t}\int_{\mid x\mid\geq 1}G(Y_{s_{-}},x)N(ds,dx),
\end{align*}
where $b,\sigma,F,G$ are some regular functions, $W$- a standard
Wiener process and $N,\tilde{N}$ - a Poisson random measure and its
compensated measure respectively. We focus on the strong
approximations, i.e. the error is measured by
$\mathbf{E}\sup_t|Y_t-\bar{Y}_t|^2$. The strong approximation is of
order $\gamma$ if $\mathbf{E}\sup_t|Y_t-\bar{Y}_t|^2\leq
\delta^{2\gamma}$, where $\delta$ is the mesh of partition of the
interval $[0,T]$. Our aim is to construct the strong approximation
for a previously fixed number $\gamma>0$. The idea is to apply the
It\^o formula to the process $Y$ many times, i.e. to the process $Y$
and then to the coefficients in its expansion. The approximation is
built of some of the coefficients which are chosen appropriately.
The main result is Theorem \ref{glowne twierdzenie} providing the
description of the approximation. This theorem is a generalization
of the results from \cite{KloPLa} for diffusion processes and
\cite{Gar} for diffusion processes with jumps generated by a
standard Poisson process. For $\gamma=\frac{1}{2}$ we obtain the
Euler approximation but we can also built approximations of higher
order. The approximation given by Theorem \ref{glowne twierdzenie}
has one limitation - in case when the L\'evy measure of a unit ball
is infinite, some ingredients are hard to simulate. This difficulty
concerns the possibility of simulating integrals with respect to the
compensated Poisson measure on unit balls. This problem hasn't
appeared in \cite{KloPLa} or \cite{Gar} since there were no jumps or
were equal to $1$ only. To overcome this problem we modify the
approximation by replacing all unit balls with $\varepsilon$- discs
which are obtained by cutting  $\varepsilon$- balls from unit balls.
This procedure causes that the error depends not only on $\delta$
but on $\varepsilon$ as well. Theorem \ref{tw glowne z epsilon}
provides the error description. It is a sum of $\delta^{2\gamma}$
and some function of $\varepsilon$ which tends to zero when
$\varepsilon\longrightarrow 0$. The speed of convergence of this
function depends on the behavior of the L\'evy measure near $0$.
Concluding, if the L\'evy measure is finite then the approximation
is given by Theorem \ref{glowne twierdzenie}, if it is not - by
Theorem \ref{tw glowne z epsilon}, but then the error depends on
$\varepsilon$ also. Note that in the first case we are able to
construct strong approximations of higher order than the Euler
scheme.

The paper is organized as follows, in Section \ref{Basic definitions
and facts} we present known facts concerning L\'evy-type stochastic
differential equation and describe the procedure of solution
expansion with the use of the It\^o formula. Section \ref{Problem
formulation} contains precise formulation of the problem which is
being successively solved in Section \ref{Construction of strong
approximation}. This section consists of three preceding lemmas
which are used in the main Theorem \ref{glowne twierdzenie}. In this
section we adopt some ideas and estimation from \cite{KloPLa} to the
present jump-diffusion settings. Section \ref{Infinite Levy measure}
is devoted to the modification of the approximation in the case
where the L\'evy measure is infinite. Section \ref{Examples}
consists of two examples of strong approximations schemes for
$\gamma=\frac{1}{2}$ and $\gamma=1$, i.e. the Euler and Milstein
schemes.

\section{Basic definitions and facts}\label{Basic definitions and facts}
Let $(\Omega,\mathcal{F}_t;t\in[0,T],P)$ be a probability space with
filtration generated by two independent processes: a standard Wiener
process $W$ and a random Poisson measure $N$. The Poisson random
measure defined on $\mathbb{R}_{+}\times(\mathbb{R}\backslash
\{0\})$ is assumed to have the intensity measure $\nu$ which is a
L\'evy measure. By $\tilde{N}$ we denote the compensated Poisson
random measure. Since we will consider stochastic integrals of
different types, the class of integrands should be specified. While
the integrals with respect to time and the Poisson measure are well
understood, the class of integrands with respect to $W$ and
$\tilde{N}$ should be made precise.
\begin{dfn}\label{def calki W}
A mapping $g_1:\Omega\times [0,T]\longrightarrow\mathbb{R}$ is
integrable with respect to $W$ if it is predictable and satisfies
the integrability condition:
$\mathbf{E}\int_{0}^{T}g_1^2(s)ds<\infty$.
\end{dfn}
\begin{dfn}\label{def calki tilde N}
Let $E$ be a subset of $\mathbb{R}$. A mapping $g_2:\Omega\times
[0,T]\times E\longrightarrow \mathbb{R}$ is integrable with respect
to $\tilde{N}$ if it is predictable and satisfies the integrability
condition:\\
$\mathbf{E}\int_{0}^{T}\int_{E}g_2^2(s,x)\nu(dx)ds<\infty$.
\end{dfn}
In these classes of integrands both integrals are square-integrable
martingales and the following {\it isometric formulas} hold:
\begin{align*}
\mathbf{E}\Big(\int_{0}^{T}g_1(s)dW_s\Big)^2&=\mathbf{E}\int_{0}^{T}g_1^2(s)ds\\
\mathbf{E}\Big(\int_{0}^{T}\int_{E}g_2(s,x)\tilde{N}(ds,dx)\Big)^2&=
\mathbf{E}\int_{0}^{T}\int_{E}g_2^2(s,x)\nu(dx)ds.
\end{align*}

\noindent Throughout all the paper we will work with a stochastic
differential equation of the form:
\begin{align}\label{rownanie_glowne}
Y_t=Y_0+\int_{0}^{t}b(Y_{s_{-}})ds+\int_{0}^{t}\sigma(Y_{s_{-}})dW_{s}&+
\int_{0}^{t}\int_{B}F(Y_{s_{-}},x)\tilde{N}(ds,dx)\nonumber\\
&+\int_{0}^{t}\int_{B^{'}}G(Y_{s_{-}},x)N(ds,dx),
\end{align}
where $t\in[0,T]$, $B=\{x: |x|<1\}$, $B^{'}=\{x: |x|\geq 1\}$. For
simplicity the initial condition is assumed to be deterministic,
i.e. $Y_0\in\mathbb{R}$. Coefficients
$b:\mathbb{R}\longrightarrow\mathbb{R}, \quad
\sigma:\mathbb{R}\longrightarrow\mathbb{R},\quad
F:\mathbb{R}\times\mathbb{R}\longrightarrow\mathbb{R},\quad
G:\mathbb{R}\times\mathbb{R}\longrightarrow\mathbb{R}$ are
measurable and satisfy the following conditions.\\
{\bf (A1)} {\bf Lipschitz condition}: there exists a constant
$K_1>0$ such that:
\begin{align*}
\mid b(y_1)-b(y_2)\mid^2&+\mid \sigma(y_1)-\sigma(y_2)\mid^2
+\int_{B}\mid
F(y_1,x)-F(y_2,x)\mid^2\nu(dx)\\[2ex]
&+\int_{B^{'}}\mid G(y_1,x)-G(y_2,x)\mid^2\nu(dx)\leq K_1\mid
y_1-y_2\mid^2 \quad \forall \ y_1,y_2\in\mathbb{R}.
\end{align*}
{\bf (A2)} {\bf Growth condition}: there exists a constant $K_2>0$
such that:
\begin{align*}
\mid b(y)\mid^2+\mid\sigma(y)\mid^2&+\int_{B}\mid
F(y,x)\mid^2\nu(dx)\\[2ex]
&+\int_{B^{'}}\mid G(y,x)\mid^2\nu(dx) \leq K_2\mid 1+y^2\mid \quad
\forall \ y\in\mathbb{R}.
\end{align*}

\begin{tw}\label{tw o istnieniu}
Under assumptions ({\bf A1}) and ({\bf A2}) there exists a unique,
adapted, c\`adl\`ag solution of (\ref{rownanie_glowne}). Moreover,
the solution satisfies:
\begin{gather}\label{drugi moment rozw}
\mathbf{E}\mid Y_t\mid^2\leq C_1(1+Y_0^2)\quad \forall \ t\in[0,T],
\end{gather}
where $C_1\geq 0$.
\end{tw}
\noindent Theorem \ref{tw o istnieniu} is a consequence of Theorem
$6.2.3$, Theorem $6.2.9$ and Corollary $6.2.4$ in \cite{Apl}, where
the Lipschitz and the growth conditions are imposed on the
coefficients $b,\sigma,F$ only and $G(\cdot,x)$ is assumed to be
continuous. The estimation $(\ref{drugi moment rozw})$ itself is a
consequence of Corollary $6.2.4$ in \cite{Apl} and the proof of
Theorem $6.2.3$, where the inequality:
\begin{gather}
\mathbf{E}\mid Y_t\mid^2\leq C(t)(1+Y_0^2)\quad \forall \ t\in[0,T],
\end{gather}
is shown for equation (\ref{rownanie_glowne}) but without the term
$\int_{0}^{t}\int_{B^{'}}G(Y_{s_{-}},x)N(ds,dx)$. Under assumptions
$\bf(A1),(A2)$ the same estimation can be obtained for
(\ref{rownanie_glowne}) with the use of similar arguments. Moreover,
$C(\cdot)$ is a continuous function and as such it is bounded on the
interval $[0,T]$ and thus $(\ref{drugi moment rozw})$ holds.

In the sequel the proposition below will be used and for the
reader's convenience we provide the proof.

\begin{prop}\label{sup drugi moment rozw }
Under assumptions ({\bf A1}) and ({\bf A2}) the solution $Y$ of
$(\ref{rownanie_glowne})$ satisfies the estimation:
\begin{gather*}
\mathbf{E}\sup_{0\leq s\leq T}\mid Y_s\mid^2\leq C_2(1+Y_0^2)
\end{gather*}
for some constant $C_2\geq 0$.
\end{prop}
{\bf Proof:} We write the solution in the form:
\begin{align*}
Y_s=Y_0&+\int_{0}^{s}b(Y_{u-})du+
\int_{0}^{s}\sigma(Y_{u-})dW_u+\int_{0}^{s}\int_{B}F(Y_{u-},x)\tilde{N}(du,dx)\\[2ex]
&+\int_{0}^{s}\int_{B^{'}}G(Y_{u-},x)\tilde{N}(du,dx)+
\int_{0}^{s}\int_{B^{'}}G(Y_{u-},x)\nu(dx)du,
\end{align*}
and thus:
\begin{align*}
Y_s^2\leq 6 \bigg(Y_0^2&+\Big(\int_{0}^{s}b(Y_{u-})du\Big)^2+
\Big(\int_{0}^{s}\sigma(Y_{u-})dW_u\Big)^2+\Big(\int_{0}^{s}\int_{B}F(Y_{u-},x)\tilde{N}(du,dx)\Big)^2\\[2ex]
&+\Big(\int_{0}^{s}\int_{B^{'}}G(Y_{u-},x)\tilde{N}(du,dx)\Big)^2+
\Big(\int_{0}^{s}\int_{B^{'}}G(Y_{u-},x)\nu(dx)du\Big)^2\bigg).
\end{align*}
Using the Doob and Schwarz inequalities as well as isometric
formulas for stochastic integrals we obtain:
\begin{align*}
\mathbf{E}\sup_{0\leq s\leq T}Y_s^2\leq 6 \bigg( Y_0^2
&+T\mathbf{E}\int_{0}^{T}b^2(Y_{u-})du+4\mathbf{E}\int_{0}^{T}\sigma^2(Y_{u-})du
+4\mathbf{E}\int_{0}^{T}\int_{B}F^2(Y_{u-},x)\nu(dx)du\\[2ex]
&+4\mathbf{E}\int_{0}^{T}\int_{B^{'}}G^2(Y_{u-},x)\nu(dx)du +
T\nu(B^{'})\mathbf{E}\int_{0}^{T}\int_{B^{'}}G^2(Y_{u-},x)\nu(dx)du
\bigg).
\end{align*}
Using assumption $\bf(A2)$ we obtain:
\begin{gather*}
\mathbf{E}\sup_{0\leq s\leq T}Y_s^2\leq 6 \bigg(
Y_0^2+K_2(T+12+T\nu(B^{'}))\int_{0}^{T}(1+\mathbf{E}Y^2_{u})du\bigg).
\end{gather*}
By (\ref{drugi moment rozw}) we have:
\begin{gather*}
\mathbf{E}\sup_{0\leq s\leq T}Y_s^2\leq 6 \bigg(
Y_0^2+K_2(T+12+T\nu(B^{'}))\int_{0}^{T}\Big(1+C_1(1+Y_0)\Big)du\bigg),
\end{gather*}
and finally we have the desired estimation:
\begin{gather*}
\mathbf{E}\sup_{0\leq s\leq T}Y_s^2\leq C_2(1+Y_0^2).
\end{gather*}
\hfill$\square$

\noindent For the process $Y$ being a solution of
(\ref{rownanie_glowne}) and for a real function $f$ of class $C^2$
we have the following form of the It\^o formula:
\begin{align}\label{wzor_Ito_glowny}
f(Y_t)=f(Y_0)&+\int_{0}^{t}f^{'}(Y_{s_-})b(Y_{s_-})ds+
\int_{0}^{t}f^{'}(Y_{s_-})\sigma(Y_{s_-})dW_s+
\frac{1}{2}\int_{0}^{t}f^{''}(Y_{s_-})\sigma^2(Y_{s_-})ds\nonumber\\[2ex]
&+\int_{0}^{t}\int_{B^{'}}\left\{f(Y_{s_-}+G(Y_{s_-},x))-f(Y_{s_-})\right\}N(ds,dx)\nonumber\\[2ex]
&+\int_{0}^{t}\int_{B}\left\{f(Y_{s_-}+F(Y_{s_-},x))-f(Y_{s_-})\right\}\tilde{N}(ds,dx)\\[2ex]
&+\int_{0}^{t}\int_{B}\left\{f(Y_{s_-}+F(Y_{s_-},x))-f(Y_{s_-})-F(Y_{s_-},x)f^{'}(Y_{s_-})\right\}\nu(dx)ds.\nonumber
\end{align}
\noindent Introducing the following operators:
\begin{align*}
&L^0f(y)&&:=f^{'}(y)b(y)+\frac{1}{2}f^{''}(y)\sigma^2(y)+\int_{B}\{f(y+F(y,x))-f(y)-F(y,x)f^{'}(y)\}\nu(dx)\\[2ex]
&L^1f(y)&&:=f^{'}(y)\sigma(y)\\[2ex]
&L^2f(y,x)&&:=f(y+F(y,x))-f(y)\\[2ex]
&L^3f(y,x)&&:=f(y+G(y,x))-f(y),
\end{align*}
we can write (\ref{wzor_Ito_glowny}) in the operator form:
\begin{align*}
f(Y_t)=f(Y_0)&+\int_{0}^{t}L^0f(Y_{s-})ds+\int_{0}^{t}L^1f(Y_{s-})dW_s\\[2ex]
&+\int_{0}^{t}\int_{B}L^2f(Y_{s-},x)\tilde{N}(ds,dx)+
\int_{0}^{t}\int_{B^{'}}L^3f(Y_{s-},x)N(ds,dx).
\end{align*}
We would like to apply the It\^o formula not only to the function
$f$, but to the coefficient functions: $L^0f$,$L^1f$,$L^2f$,$L^3f$
or in general to any function which is smooth enough as well. Since
functions $L^2f$ and $L^3f$ depend on two arguments $(x,y)$, we
admit the following rules of acting operators on the multiargument
real function $g(y,x_1,x_2,...,x_l)$:
\begin{align*}
&L^0g(y,x_1,...,x_l):=\frac{\partial}{\partial y}g(y,x_1,...,x_l)b(y)+\frac{1}{2}\frac{\partial^2}{\partial y^2}g(y,x_1,...,x_l)\sigma^2(y)\\[2ex]
&\phantom{L^0}+\int\limits_{B}\Big\{g(y+F(y,x),x_1,...,x_l)-g(y,x_1,...,x_l)-F(y,x)\frac{\partial}{\partial y}g(y,x_1,...,x_l)\Big\}\nu(dx)\\[2ex]
&L^1g(y,x_1,...,x_l):=\frac{\partial}{\partial y}g(y,x_1,...,x_l)\sigma(y)\\[2ex]
&L^2g(y,x_1,...,x_l,x_{l+1}):=g(y+F(y,x_{l+1}),x_1,...,x_l)-g(y,x_1,...,x_l)\\[2ex]
&L^3g(y,x_1,...,x_l,x_{l+1}):=g(y+G(y,x_{l+1}),x_1,...,x_l)-g(y,x_1,...,x_l).
\end{align*}
\noindent To describe the higher order It\^o expansion of $f$ we
will use the notion of {\it multiindices} and {\it multiple
stochastic integrals}. A multiindex
$\alpha=(\alpha_1,\alpha_2,...,\alpha_{l(\alpha))}$ is a finite
sequence of elements such that $\alpha_i\in\{0,1,2,3\}$ for
$i=1,2,...,l(\alpha)$. The number of all elements equal to
\begin{enumerate}[a)]
\item $0$ will be denoted by $s(\alpha)$,
\item $1$ will be denoted by $w(\alpha)$,
\item $2$ will be denoted by $\tilde{n}(\alpha)$,
\item $3$ will be denoted by $n(\alpha)$.
\end{enumerate}
 The length $l(\alpha)$ of
$\alpha$ is thus given as
$l(\alpha)=s(\alpha)+w(\alpha)+\tilde{n}(\alpha)+n(\alpha)$. For the
sake of convenience we also define
$k(\alpha):=\tilde{n}(\alpha)+n(\alpha)$. For technical reasons we
also consider the empty index denoted by $v$ with length $0$, i.e.
$l(v)=0$. For a given multiindex
$\alpha=(\alpha_1,\alpha_2,...,\alpha_{l(\alpha)})$ let us define:
\begin{align*}
\alpha-&=(\alpha_1,\alpha_2,...,\alpha_{l(\alpha)-1})\\[2ex]
-\alpha&=(\alpha_2,...,\alpha_{l(\alpha)}) .
\end{align*}
\begin{dfn}
A set of multiindices $\mathcal{A}$ is called a {\bf hierarchical
set} if $\forall \alpha\in\mathcal{A}:$
\begin{gather*}
 l(\alpha)<\infty \quad \text{and} \quad \alpha\in\mathcal{A}\backslash \{v\} \Longrightarrow -\alpha\in\mathcal{A}.
\end{gather*}
A set of multiindices $\mathcal{B}(\mathcal{A})$, where
$\mathcal{A}$ is a hierarchical set, is called a {\bf remainder set}
of $\mathcal{A}$ if \ $\forall \alpha\in\mathcal{B}(\mathcal{A})$
\begin{gather*}
\alpha\notin\mathcal{A}\quad \text{and} \quad -\alpha\in\mathcal{A}.
\end{gather*}
\end{dfn}
\noindent Assume that $g(s,x_1,x_2,...,x_l)$ is a regular stochastic
process, i.e. such that all the stochastic integrals written below
exist in the sense of Definitions \ref{def calki W} and \ref{def
calki tilde N}. Let $\rho$ and $\tau$ be fixed points in the
interval $[0,T]$ s.t. $\rho\leq\tau$. A multiple stochastic integral
on the interval $[\rho,\tau]$ with respect to any multiindex
$\alpha$ s.t. $k(\alpha)\leq l$ is defined by the induction
procedure. First, we define the integral with respect to the empty
index:
\begin{gather*}
I_{v}[g]_{\rho}^{\tau}(x_1,...,x_l)=g(\tau,x_1,...,x_l).
\end{gather*}
Now, assume that $I_{\alpha-}[g]_{\rho}^{\tau}(x_1,x_2,...,x_k)$
depends on $k$ parameters, where $0\leq k\leq l$. Then we define the
multiple integral as follows:
\begin{enumerate}[1)]
\item if $\alpha_{l(\alpha)}=0$ then
\begin{gather*}
I_{\alpha}[g]_{\rho}^{\tau}(x_1,...,x_{k})=\int_{\rho}^{\tau}I_{\alpha-}[g]_{\rho}^{s-}(x_1,...,x_{k})ds,
\end{gather*}
\item if $\alpha_{l(\alpha)}=1$ then
\begin{gather*}
I_{\alpha}[g]_{\rho}^{\tau}(x_1,...,x_{k})=\int_{\rho}^{\tau}I_{\alpha-}[g]_{\rho}^{s-}(x_1,...,x_{k})dW_s,
\end{gather*}
\item if $\alpha_{l(\alpha)}=2$ and $k\geq 1$ then
\begin{gather*}
I_{\alpha}[g]_{\rho}^{\tau}(x_1,...,x_{k-1})=\int_{\rho}^{\tau}\int_{B}I_{\alpha-}[g]_{\rho}^{s-}(x_1,...,x_{k})\tilde{N}(ds,dx_{k}),
\end{gather*}
\item if $\alpha_{l(\alpha)}=3$ and $k\geq 1$ then
\begin{gather*}
I_{\alpha}[g]_{\rho}^{\tau}(x_1,...,x_{k-1})=\int_{\rho}^{\tau}\int_{B^{'}}I_{\alpha-}[g]_{\rho}^{s-}(x_1,...,x_{k})N(ds,dx_{k}).
\end{gather*}
\end{enumerate}
Let us notice that it follows from the description above that
$I_{\alpha}[g]$ depends on $l-k(\alpha)$ parameters, i.e.
$I_{\alpha}[g]_{\rho}^{\tau}=I_{\alpha}[g]_{\rho}^{\tau}(x_1,x_2,...,x_{l-k(\alpha)})$.\\
\\
\noindent {\bf Example} Let $g=g(s,x_1,x_2,x_3)$. Then:
\begin{align*}
&I_{(1)}[g]_{\rho}^{\tau}(x_1,x_2,x_3)=\int_{\rho}^{\tau}g(s-,x_1,x_2,x_3)dW_s,\\[2ex]
&I_{(213)}[g]_{0}^{t}(x_1)=\int\limits_{0}^{t}\int\limits_{B^{'}} \
\ \int\limits_{0}^{\ \ s_1-} \ \ \int\limits_{0}^{ \
 \ s_2-}\int\limits_{B}g(s_3-,x_1,x_2,x_3) \ \tilde{N}(ds_3,dx_3) \
dW_{s_2} \ N(ds_1,dx_2).
\end{align*}

\noindent The processes which serve as integrands in multiple
integrals in the expansion of $f(Y)$ will be obtained with the use
of {\it coefficient functions} $f_{\alpha}$, where $\alpha$ is a
multiindex. We define the coefficient function with respect to any
multiindex $\alpha$  by the induction procedure:
\begin{align*}
f_v(y)&=f(y),\\[2ex]
f_{\alpha}(y,x_1,...,x_{k(\alpha)})&=L^{\alpha_{1}}\Big[f_{-\alpha}(y,x_1,...,x_{k(-\alpha)})\Big](y,x_1,...,x_{k(\alpha)}).
\end{align*}

\noindent {\bf Example} For a given function $f=f(y)$ we get:
\begin{align*}
f_{(10)}(y)&=L^1L^0f,\\[2ex]
f_{(2013)}(y,x_1,x_2)&=L^2L^0L^1L^3f.
\end{align*}
For simplicity we omit here the dependence on arguments on the right
hand side.\\

\noindent Notice, that the coefficient function
$f_\alpha=f_\alpha(y,x_1,...,x_{k(\alpha)})$ depends on $k(\alpha)$
parameters, i.e. on $x_1,x_2,...,x_{k(\alpha)}$. However, the
multiple integral
$I_{\alpha}[f_{\alpha}]_{\rho}^{\tau}=I_{\alpha}[f_\alpha(y,x_1,...,x_{k(\alpha)})]_{\rho}^{\tau}$
does not depend on any parameter.\\

\noindent We have the following analogue of Theorem 5.5.1 in
\cite{KloPLa} which is also called the It\^o - Taylor expansion. It
is a consequence of the It\^o formula and definitions of the
hierarchical and remainder sets.
\begin{tw}
For any hierarchical set $\mathcal{A}$ and a smooth function $f$ we
have the following representation:
\begin{gather}\label{wzor_Ito_hierarchiczny}
f(Y_\tau)=\sum_{\alpha\in\mathcal{A}}I_{\alpha}[f_{\alpha}(Y_\rho,x_1,...,x_{k(\alpha)})]_{\rho}^{\tau}+
\sum_{\alpha\in\mathcal{B(\mathcal{A})}}I_{\alpha}[f_{\alpha}(Y_{\bullet-},x_1,...,x_{k(\alpha)})]_{\rho}^{\tau},
\end{gather}
\noindent assuming that all the integrals above exist.
\end{tw}
Notice that the first sum in (\ref{wzor_Ito_hierarchiczny}) consists
of all integrals for which the integrands do not depend on time
while the second sum contains all integrals with the integrands
dependent on time. \noindent Since we are interested in the
approximation of the process $Y$ itself, to the end of the paper we
will consider the identity function only, i.e. $f(y)=y$.\\

\noindent In the sequel we use two auxiliary lemmas.
\begin{lem}[The Gronwall lemma]\label{Gronwal lemma}
Let $g,h:[0,T]\longrightarrow \mathbb{R}$ be integrable and satisfy:
\begin{gather*}
0\leq g(t)\leq h(t)+L\int_{0}^{t}g(s)ds
\end{gather*}
for $t\in[0,T]$ and $L>0$. Then:
\begin{gather*}
g(t)\leq h(t)+L\int_{0}^{t}e^{L(t-s)}h(s)ds
\end{gather*}
for $t\in[0,T].$
\end{lem}
\begin{lem}\label{cadlag lemma}
Let $g$ be a c\`adl\`ag function on the interval $[0,T]$. Then for
any $(\rho,\tau]\subseteq [0,T]$ we have:
\begin{gather*}
\sup_{s\in(\rho,\tau]} g(s-)\leq\sup_{s\in(\rho,\tau]} g(s).
\end{gather*}
\end{lem}
{\bf Proof:} Let $(s_n)_{n=1,2,...}$ be a sequence such that
$s_n\in(\rho,\tau]$ for $n=1,2,...$ satisfying
$g(s_n-)\longrightarrow\sup_{s\in(\rho,\tau]}g(s-):=K$. Since $g$ is
c\`adl\`ag, for any $\varepsilon>0$ there exits a sequence
$(s^{\varepsilon}_{n})_{n=1,2,...}$ such that
$s_n^{\varepsilon}\in(\rho,\tau]$ for $n=1,2,...$ and satisfies:
\begin{gather*}
g(s^{\varepsilon}_{n})\geq g(s_n-)-\varepsilon \quad \text{for} \
n=1,2,...
\end{gather*}
and thus:
\begin{gather*}
\underset{n\longrightarrow
\infty}{\overline{\lim}}g(s_n^{\varepsilon})\geq K-\varepsilon.
\end{gather*}
Letting $\varepsilon\longrightarrow 0$ we obtain
$sup_{s\in(\rho,\tau]}g(s)\geq K.$ \hfill$\square$

\section{Problem formulation}\label{Problem formulation}
Our approximation of the process $Y$, which is the solution of
(\ref{rownanie_glowne}), will be based on a fixed partition
\begin{gather*}
0=\tau_0<\tau_1<...<\tau_n=T
\end{gather*}
of the interval $[0,T]$. For the sake of simplicity all the
partition points are assumed to be non-random. The diameter of this
partition is assumed to be smaller than $\delta$, i.e.\\
$\max_{i=0,1,...,n-1}(\tau_{i+1}-\tau_i)<\delta$. The approximation
denoted by $Y^{\delta}$ is obtained from the first sum of multiple
integrals in the It\^o-Taylor expansion
(\ref{wzor_Ito_hierarchiczny}). The procedure can be described as
follows. Starting from the known value $Y^{\delta}_{0}$, which can
be equal to $Y_0$, we calculate the value $Y^{\delta}_t$ for $t\in
(0,\tau_1]$ using the first sum in (\ref{wzor_Ito_hierarchiczny}).
Using value $Y^{\delta}_{\tau_1}$ we repeat the procedure for
$t\in(\tau_1, \tau_2]$ and so on. Denoting $n_t=\max\{k : \tau_k\leq
t\}$ we define process $Y^{\delta}$ as:
\begin{gather}\label{wzor_postac_aprox.}
Y^{\delta}_t=\sum_{\alpha\in\mathcal{A}}I_{\alpha}[f_{\alpha}(Y^\delta_{\tau_{n_t}},x_1,...,x_{k(\alpha)})]_{\tau_{n_t}}^{t}.
\end{gather}
The motivation for the form of the approximation is justified by the
possibility of practical calculation multiple integrals for which
integrands does not depend on time (at least for low order
integrals). In fact, in the case of integrals with respect to the
compensated Poisson measure additional difficulty occurs which is
related to the property of L\'evy measure. It is discussed in
Section \ref{Infinite Levy measure}. We focus on the problem of
finding a strong approximation of order $\gamma>0$, i.e. such that
\begin{gather}\label{def_bledu}
\mathbb{E}\sup_{t\in[0,T]}\mid Y_t-Y^\delta_t\mid^2\leq
C\delta^{2\gamma}
\end{gather}
for some constant $C>0$. The rate of convergence $\gamma$ is fixed
and in practical application it is the multiplicity of
$\frac{1}{2}$, i.e. $\gamma=\frac{1}{2},1,\frac{3}{2},...$.
\\

\noindent Thus our goal can be summarized as follows: for a fixed
$\gamma>0$ find a hierarchical set $\mathcal{A}$ such that the
approximation $Y^\delta$ defined by (\ref{wzor_postac_aprox.})
satisfies (\ref{def_bledu}).

\section{Construction of the strong approximation}\label{Construction of strong approximation}
Before formulating the main theorem let us introduce the following
notation. For any multiindex $\alpha$ s.t. $k(\alpha)>0$  we denote
by $\beta(\alpha)$ a multiindex which is obtained from $\alpha$ by
deleting all the coordinates equal to $0$ or $1$. Then the sets
$B_i^{\alpha}$ for $i=1,2,...,k(\alpha)$ are defined as follows
\[B_i^{\alpha}:= \
\begin{cases}
B  \ \text{if} \ \beta(\alpha)_{k(\alpha)+1-i}=2
\\
B^{'} \ \text{if} \ \beta(\alpha)_{k(\alpha)+1-i}=3.
\end{cases}
\]
Recall that $B$ is a unit ball and $B^{'}$ its complement. The
following result is a generalization of Theorem $10.6.3$ in
\cite{KloPLa} and Theorem $7$ in \cite{Gar}.
\begin{tw}\label{glowne twierdzenie}
Let us assume that coefficients in equation (\ref{rownanie_glowne})
satisfy conditions {\bf(A1)},{\bf(A2)}. Let $Y^{\delta}$ be the
approximation of the form (\ref{wzor_postac_aprox.}), for the
solution $Y$ of (\ref{rownanie_glowne}), constructed with the use of
the hierarchical set $\mathcal{A}_\gamma$, where:
\begin{gather}\label{hierarchical set}
\mathcal{A}_{\gamma}:=\left\{ \ l(\alpha)+s(\alpha)\leq 2\gamma
\quad \text{or} \quad l(\alpha)=s(\alpha)=\gamma+\frac{1}{2} \
\right\}.
\end{gather}
Moreover, assume that
coefficient functions $f_{\alpha}$ satisfy:\\
{\bf(A3)} for any $\alpha\in\mathcal{A}_{\gamma}$ holds:
\begin{align*}
\int\limits_{B_1^{\alpha}}\int\limits_{B_2^{\alpha}}...\int\limits_{B_{k(\alpha)}^{\alpha}}\mid
f_{\alpha}(y_1,x_1,...,x_{k(\alpha)})-f_{\alpha}(y_2,x_1,...,x_{k(\alpha)})\mid^2\nu(dx_{k(\alpha)})...\nu(dx_1)\\[2ex]
\leq K_{\alpha} \mid y_1-y_2\mid^2,
\end{align*}
{\bf(A4)} for any
$\alpha\in\mathcal{A}_{\gamma}\cup\mathcal{B}(\mathcal{A}_{\gamma})$
holds:
\begin{align*}
\int_{B_1^{\alpha}}\int_{B_2^{\alpha}}...\int_{B_{k(\alpha)}^{\alpha}}\mid
f_{\alpha}(y,x_1,x_2,...,x_{k(\alpha)})\mid^2\nu(dx_{k(\alpha)})...\nu(dx_1)&\leq
L_{\alpha} ( 1+ y^2 ),
\end{align*}
 where
$K_\alpha$, $L_{\alpha}$ are some constants.\\
Then for $\delta\in(0,1)$ the inequality:
\begin{gather*}
\mathbf{E}\sup_{s\in[0,T]}\mid Y_s-Y^{\delta}_s\mid^2\leq
E_1(\gamma,T)\mid Y_0-Y^{\delta}_0
\mid^2+E_2(\gamma,T,Y_0)\delta^{2\gamma}
\end{gather*}
holds.
\end{tw}

\noindent The proof is presented at the end of this section. First
we present three auxiliary lemmas and a proposition.
\begin{lem}\label{lemat1}
Let $\rho,\tau$ be two fixed points in the interval $[0,T]$ s. t.
$\rho<\tau$, $\tau-\rho<\delta$. If all the integrals below exist
then we have:

\begin{align}\label{L11}
\mathbf{E}\sup_{s\in(\rho,\tau]}\Big\{\int_{\rho}^{s}g(u,x_1,x_2,...,x_l)du\Big\}^2&\leq
\delta^2 \mathbf{E}\Big\{\sup_{u\in(\rho,\tau]}g^2(u,x_1,x_2,...,x_l)\Big\},\\[2ex]\label{L12}
\mathbf{E}\sup_{s\in(\rho,\tau]}\Big\{\int_{\rho}^{s}g(u,x_1,x_2,...,x_l)du\Big\}^2&\leq\delta
\int_{\rho}^{\tau}\mathbf{E}\big\{g^2(u,x_1,x_2,...,x_l)\big\}du,
\end{align}

\vskip5pt

\begin{align}\label{L13}
\mathbf{E}\sup_{s\in(\rho,\tau]}\Big\{\int_{\rho}^{s}g(u,x_1,x_2,...,x_l)dW_u\Big\}^2&\leq
4\delta \ \mathbf{E}\Big\{\sup_{u\in(\rho,\tau]}g^2(u,x_1,x_2,...,x_l)\Big\},\\[2ex]\label{L14}
\mathbf{E}\sup_{s\in(\rho,\tau]}\Big\{\int_{\rho}^{s}g(u,x_1,x_2,...,x_l)dW_u\Big\}^2&\leq
4 \
\int_{\rho}^{\tau}\mathbf{E}\big\{g^2(u,x_1,x_2,...,x_l)\big\}du,
\end{align}

\vskip5pt

\begin{align}\label{L15}
\mathbf{E}\sup_{s\in(\rho,\tau]}\Big\{\int_{\rho}^{s}\int_{B}g(u,x_1,...,x_l)\tilde{N}(du,dx_l)\Big\}^2&\leq
4\delta \ \int_{B}\mathbf{E}\Big\{\sup_{u\in(\rho,\tau]}g^2(u,x_1,...,x_l)\Big\}\nu(dx_l),\\[2ex]\label{L16}
\mathbf{E}\sup_{s\in(\rho,\tau]}\Big\{\int_{\rho}^{s}\int_{B}g(u,x_1,...,x_l)\tilde{N}(du,dx_l)\Big\}^2&\leq
4 \
\int_{\rho}^{\tau}\int_{B}\mathbf{E}\big\{g^2(u,x_1,x_2,...,x_l)\big\}\nu(dx_l)du,
\end{align}

\vskip5pt

\begin{align}\label{L17}
\mathbf{E}\sup_{s\in(\rho,\tau]}\Big\{\int_{\rho}^{s}\int_{B^{'}}g(u,x_1,...,x_l)&N(du,dx_l)\Big\}^2\leq\nonumber\\[2ex]
&2\delta(4+\delta\nu(B^{'})) \ \int_{B^{'}}\mathbf{E}\Big\{\sup_{u\in(\rho,\tau]}g^2(u,x_1,...,x_l)\Big\}\nu(dx_l),\\[2ex]\label{L18}
\mathbf{E}\sup_{s\in(\rho,\tau]}\Big\{\int_{\rho}^{s}\int_{B^{'}}g(u,x_1,...,x_l)N&(du,dx_l)\Big\}^2\leq\nonumber\\[2ex]
&2(4+\delta\nu(B^{'})) \
\int_{\rho}^{\tau}\int_{B^{'}}\mathbf{E}\big\{g^2(u,x_1,...,x_l)\big\}\nu(dx_l)du.
\end{align}

\end{lem}
\noindent Note, that due to Lemma \ref{cadlag lemma}, the lemma
above remains true if we replace the upper limit "$s$" in the left
hand
side integrals with $"s-"$.\\

\noindent {\bf Proof:} All these inequalities are proved with the
use of the Schwarz and Doob inequalities, the isometric formula for
stochastic integrals and Fubini's theorem.\\
\noindent (\ref{L12})
\begin{align*}
\mathbf{E}\sup_{s\in(\rho,\tau]}\Big\{\int_{\rho}^{s}g(u,x_1,x_2,...,x_l)du\Big\}^2&\leq\mathbf{E}\sup_{s\in(\rho,\tau]}\delta\Big\{\int_{\rho}^{s}g^2(u,x_1,x_2,...,x_l)du\Big\}\\[2ex]
\leq\delta\mathbf{E}\Big\{\int_{\rho}^{\tau}g^2(u,x_1,x_2,...,x_l)du\Big\}&=\delta
\ \int_{\rho}^{\tau}\mathbf{E}\big\{g^2(u,x_1,x_2,...,x_l)\big\}du
\end{align*}
(\ref{L11})
\begin{align*}
\mathbf{E}\sup_{s\in(\rho,\tau]}\Big\{\int_{\rho}^{s}g(u,x_1,x_2,...,x_l)du\Big\}^2&\leq\delta
\
\int_{\rho}^{\tau}\mathbf{E}\big\{g^2(u,x_1,x_2,...,x_l)\big\}du\\[2ex]
&\leq\delta^2 \
\mathbf{E}\Big\{\sup_{u\in(\rho,\tau]}g^2(u,x_1,x_2,...,x_l)\Big\}
\end{align*}
(\ref{L14})
\begin{align*}
\mathbf{E}\sup_{s\in(\rho,\tau]}\Big\{\int_{\rho}^{s}g(u,x_1,x_2,...,x_l)dW_u\Big\}^2&\leq
4\sup_{s\in(\rho,\tau]}\mathbf{E}\Big\{\int_{\rho}^{s}g(u,x_1,x_2,...,x_l)dW_u\Big\}^2\\[2ex]
=4\sup_{s\in(\rho,\tau]}\mathbf{E}\Big\{\int_{\rho}^{s}g^2(u,x_1,x_2,...,x_l)du\Big\}&\leq
4\mathbf{E}\Big\{\int_{\rho}^{\tau}g^2(u,x_1,x_2,...,x_l)du\Big\}\\[2ex]
&=4\int_{\rho}^{\tau}\mathbf{E}\Big\{g^2(u,x_1,x_2,...,x_l)\Big\}du
\end{align*}
(\ref{L13})
\begin{align*}
\mathbf{E}\sup_{s\in(\rho,\tau]}\Big\{\int_{\rho}^{s}g(u,x_1,x_2,...,x_l)dW_u\Big\}^2&\leq
4\mathbf{E}\Big\{\int_{\rho}^{\tau}g^2(u,x_1,x_2,...,x_l)du\Big\}\\[2ex]
&\leq4\delta \
\mathbf{E}\Big\{\sup_{u\in(\rho,\tau]}g^2(u,x_1,x_2,...,x_l)\Big\}
\end{align*}
(\ref{L16})
\begin{align*}
\mathbf{E}\sup_{s\in(\rho,\tau]}\Big\{\int_{\rho}^{s}\int_{B}g(u,x_1,...,x_l)\tilde{N}(du,dx_l)\Big\}^2&\leq
4\sup_{s\in(\rho,\tau]}\mathbf{E}\Big\{\int_{\rho}^{s}\int_{B}g(u,x_1,...,x_l)\tilde{N}(du,dx_l)\Big\}^2\\[2ex]
=4\sup_{s\in(\rho,\tau]}\mathbf{E}\Big\{\int_{\rho}^{s}\int_{B}g^2(u,x_1,...,x_l)\nu(dx_l)du\Big\}
&\leq4 \
\int_{\rho}^{\tau}\int_{B}\mathbf{E}\big\{g^2(u,x_1,x_2,...,x_l)\big\}\nu(dx_l)du
\end{align*}
(\ref{L15})
\begin{align*}
\mathbf{E}\sup_{s\in(\rho,\tau]}\Big\{\int_{\rho}^{s}\int_{B}g(u,x_1,...,x_l)\tilde{N}(du,dx_l)\Big\}^2&\leq4
\
\mathbf{E}\Big\{\int_{\rho}^{\tau}\int_{B}g^2(u,x_1,x_2,...,x_l)\nu(dx_l)du\Big\}\\[2ex]
\leq4 \
\mathbf{E}\Big\{\delta\sup_{u\in(\rho,\tau]}\int_{B}g^2(u,x_1,x_2,...,x_l)\nu(dx_l)\Big\}&\leq
4\delta\mathbf{E}\Big\{\int_{B}\sup_{u\in(\rho,\tau]}g^2(u,x_1,x_2,...,x_l)\nu(dx_l)\Big\}\\[2ex]
&=4\delta \
\int_{B}\mathbf{E}\Big\{\sup_{u\in(\rho,\tau]}g^2(u,x_1,...,x_l)\Big\}\nu(dx_l)
\end{align*}
(\ref{L18})
\begin{align}\label{nier_zamiana_zwykl_na_skomp}
&\mathbf{E}\sup_{s\in(\rho,\tau]}\Big\{\int_{\rho}^{s}\int_{B^{'}}g(u,x_1,...,x_l)\tilde{N}(du,dx_l)\Big\}^2
\nonumber\\[2ex]
&=\mathbf{E}\sup_{s\in(\rho,\tau]}\Big\{
\int_{\rho}^{s}\int_{B^{'}}g(u,x_1,...,x_l)N(du,dx_l)
+\int_{\rho}^{s}\int_{B^{'}}g(u,x_1,...,x_l)\nu(dx_l)du
\Big\}^2\nonumber\\[2ex]
&\leq 2\Big\{\mathbf{E}\sup_{s\in(\rho,\tau]}\Big\{
\int\limits_{\rho}^{s}\int\limits_{B^{'}}g(u,x_1,...,x_l)\tilde{N}(du,dx_l)\Big\}^2
+\mathbf{E}\sup_{s\in(\rho,\tau]}\Big\{\int\limits_{\rho}^{s}\int\limits_{B^{'}}g(u,x_1,...,x_l)\nu(dx_l)du
\Big\}^2\Big\}
\end{align}
The first component is bounded by analogous expression as in
(\ref{L16}). For the second we have the following inequalities:
\begin{align*}
\mathbf{E}\sup_{s\in(\rho,\tau]}\Big\{\int_{\rho}^{s}\int_{B^{'}}&g(u,x_1,...,x_l)\nu(dx_l)du
\Big\}^2\\[2ex]
&\leq
\mathbf{E}\sup_{s\in(\rho,\tau]}\Big\{\int_{\rho}^{s}\int_{B^{'}}1\nu(dx_l)du\cdot\int_{\rho}^{s}\int_{B^{'}}g^2(u,x_1,...,x_l)\nu(dx_l)du
\Big\}^2\\[2ex]
&\leq\delta\nu(B^{'})\int_{\rho}^{\tau}\int_{B^{'}}\mathbf{E}\big\{g^2(u,x_1,...,x_l)\big\}\nu(dx_l)du.
\end{align*}
As a consequence we obtain:
\begin{gather*}
\mathbf{E}\sup_{s\in(\rho,\tau]}\Big\{\int\limits_{\rho}^{s}\int\limits_{B^{'}}g(u,x_1,...,x_l)N(du,dx_l)\Big\}^2\leq
\ 2(4+\delta\nu(B^{'})) \
\int\limits_{\rho}^{\tau}\int\limits_{B^{'}}\mathbf{E}\big\{g^2(u,x_1,...,x_l)\big\}\nu(dx_l)du.
\end{gather*}
(\ref{L17})\\
\noindent For the second term in
$(\ref{nier_zamiana_zwykl_na_skomp})$ we have the following
inequalities:
\begin{align*}
\mathbf{E}\sup_{s\in(\rho,\tau]}\Big\{\int_{\rho}^{s}\int_{B^{'}}g(u,x_1,...,x_l)\nu(dx_l)du
\Big\}^2&\leq\delta\nu(B^{'})\mathbf{E}\Big\{\int_{\rho}^{\tau}\int_{B^{'}}g^2(u,x_1,...,x_l)\nu(dx_l)du\Big\}\\[2ex]
\leq\delta^2\nu(B^{'})\mathbf{E}\Big\{\sup_{u\in(\rho,\tau]}\int_{B^{'}}g^2(u,x_1,...,x_l)\nu(dx_l)\Big\}
&\leq\delta^2\nu(B^{'})\mathbf{E}\Big\{\int_{B^{'}}\sup_{u\in(\rho,\tau]}g^2(u,x_1,...,x_l)\nu(dx_l)\Big\}\\[2ex]
&=\delta^2\nu(B^{'})\int_{B^{'}}\mathbf{E}\Big\{\sup_{u\in(\rho,\tau]}g^2(u,x_1,...,x_l)\Big\}\nu(dx_l).
\end{align*}
Taking into account (\ref{nier_zamiana_zwykl_na_skomp}), the
inequality above and $(\ref{L15})$ we obtain:
\begin{gather*}
\mathbf{E}\sup_{s\in(\rho,\tau]}\Big\{\int\limits_{\rho}^{s}\int\limits_{B^{'}}g(u,x_1,...,x_l)N(du,dx_l)\Big\}^2\leq
 2\delta(4+\delta\nu(B^{'})) \
\int\limits_{B^{'}}\mathbf{E}\Big\{\sup_{u\in(\rho,\tau]}g^2(u,x_1,...,x_l)\Big\}\nu(dx_l).
\end{gather*}
$\hfill\square$

\begin{lem}\label{lemat2}
Let $\rho,\tau$ be two fixed points in the interval $[0,T]$ s.t.
$\rho<\tau$, $\tau-\rho<\delta$ and $\alpha\neq v$ be a fixed
multiindex. If all the integrals below exist for the process
$g=g(u,x_1,...,x_l)$, where $l\geq k(\alpha)$, then we have:
\begin{align}
\mathbf{E}\{\sup_{s\in(\rho,\tau]}I^2_{\alpha}&[g]_{\rho}^{s}(x_1,...,x_{l-k(\alpha)})\}
\leq
\delta^{l(\alpha)+s(\alpha)-1}4^{w(\alpha)+\tilde{n}(\alpha)}\Big\{2(4+\delta\nu(B^{'}))\Big\}^{n(\alpha)}\cdot\nonumber\\[2ex]
\label{wz_lem2}
&\cdot\int_{\rho}^{\tau}\mathbf{E}\int\limits_{B_1^{\alpha}}\int\limits_{B_2^{\alpha}}...\int\limits_{B_{k(\alpha)}^{\alpha}}g^2(u,x_1,x_2,...x_{l})\nu(dx_{l})
\ \nu(dx_{l-1}) \ ... \ \nu(dx_{l-k(\alpha)+1}) du.
\end{align}
\end{lem}
{\bf Proof:} We will apply the induction procedure with respect to
the length of $\alpha$. If $l(\alpha)=1$ then (\ref{wz_lem2})
follows from inequalities (\ref{L12}), (\ref{L14}), (\ref{L16}),
(\ref{L18}) in Lemma \ref{lemat1} applied to $\alpha=0$, $\alpha=1$,
$\alpha=2$,
$\alpha=3$ respectively.\\
\noindent Now assume that $(\ref{wz_lem2})$ is true for $\alpha-$
and let us show that it is also true for $\alpha$. We will consider
several
cases.\\
\noindent $a) \ \alpha_{l(\alpha)}=0$; In this case
$k(\alpha-)=k(\alpha)$ and $B_i^{\alpha-}=B_i^{\alpha}$ for
$i=1,2,...,k(\alpha)$. By (\ref{L11}), Lemma \ref{cadlag lemma} and
the inductive assumption we have:
\begin{align*}
&\mathbf{E}\{\sup_{s\in(\rho,\tau]}I^2_{\alpha}[g]_{\rho}^{s}(x_1,...,x_{l-k(\alpha)})\}=
\mathbf{E}\Big\{\sup_{s\in(\rho,\tau]}\int_{\rho}^{s}I_{\alpha-}[g]_{\rho}^{u-}(x_1,...,x_{l-k(\alpha-)})du\Big\}^2\\[2ex]
&\leq\delta^2\mathbf{E}\Big\{\sup_{u\in(\rho,\tau]}I^2_{\alpha-}[g]_{\rho}^{u-}(x_1,...,x_{l-k(\alpha-)})\Big\}\\[2ex]
&\leq\delta^2
\delta^{l(\alpha-)+s(\alpha-)-1}4^{w(\alpha-)+\tilde{n}(\alpha-)}\Big\{2(4+\delta\nu(B^{'}))\Big\}^{n(\alpha-)}\cdot\\[2ex]
&\cdot\int_{\rho}^{\tau}\mathbf{E}\int\limits_{B_1^{\alpha-}}\int\limits_{B_2^{\alpha-}}...\int\limits_{B_{k(\alpha-)}^{\alpha-}}g^2(u,x_1,x_2,...,x_{l})\nu(dx_{l})
\ \nu(dx_{l-1}) \ ... \ \nu(dx_{l-k(\alpha-)+1})
du\\[2ex]
\end{align*}
\begin{align*}
&=\delta^{l(\alpha)+s(\alpha)-1}4^{w(\alpha)+\tilde{n}(\alpha)}\Big\{2(4+\delta\nu(B^{'}))\Big\}^{n(\alpha)}\cdot\\[2ex]
&\cdot\int_{\rho}^{\tau}\mathbf{E}\int\limits_{B_1^{\alpha}}\int\limits_{B_2^{\alpha}}...\int\limits_{B_{k(\alpha)}^{\alpha}}g^2(u,x_1,x_2,...,x_{l})\nu(dx_{l})
\ \nu(dx_{l-1}) \ ... \ \nu(dx_{l-k(\alpha)+1}) du.
\end{align*}

\noindent $b) \ \alpha_{l(\alpha)}=1$; In this case
$k(\alpha-)=k(\alpha)$ and $B_i^{\alpha-}=B_i^{\alpha}$ for
$i=1,2,...,k(\alpha)$. By (\ref{L13}), Lemma \ref{cadlag lemma} and
the inductive assumption we have:
\begin{align*}
&\mathbf{E}\{\sup_{s\in(\rho,\tau]}I^2_{\alpha}[g]_{\rho}^{s}\}(x_1,...,x_{l-k(\alpha)})=
\mathbf{E}\Big\{\sup_{s\in(\rho,\tau]}\int_{\rho}^{s}I_{\alpha-}[g]_{\rho}^{u-}(x_1,...,x_{l-k(\alpha-)})
dW_u\Big\}^2\\[2ex]
&\leq4\delta\mathbf{E}\Big\{\sup_{u\in(\rho,\tau]}I^2_{\alpha-}[g]_{\rho}^{u-}(x_1,...,x_{l-k(\alpha-)})\Big\}\\[2ex]
&\leq4\delta \
\delta^{l(\alpha-)+s(\alpha-)-1}4^{w(\alpha-)+\tilde{n}(\alpha-)}\Big\{2(4+\delta\nu(B^{'}))\Big\}^{n(\alpha-)}\cdot\\[2ex]
&\cdot\int_{\rho}^{\tau}\mathbf{E}\int\limits_{B_1^{\alpha-}}\int\limits_{B_2^{\alpha-}}...\int\limits_{B_{k(\alpha-)}^{\alpha-}}g^2(u,x_1,x_2,...x_{l})\nu(dx_{l})
\ \nu(dx_{l-1}) \ ... \ \nu(dx_{l-k(\alpha-)+1})
du\\[2ex]
&=\delta^{l(\alpha)+s(\alpha)-1}4^{w(\alpha)+\tilde{n}(\alpha)}\Big\{2(4+\delta\nu(B^{'}))\Big\}^{n(\alpha)}\cdot\\[2ex]
&\cdot\int_{\rho}^{\tau}\mathbf{E}\int\limits_{B_1^{\alpha}}\int\limits_{B_2^{\alpha}}...\int\limits_{B_{k(\alpha)}^{\alpha}}g^2(u,x_1,x_2,...x_{l})\nu(dx_{l})
\ \nu(dx_{l-1}) \ ... \ \nu(dx_{l-k(\alpha)+1}) du.
\end{align*}
\noindent $c) \ \alpha_{l(\alpha)}=2$; By (\ref{L15}), Lemma
\ref{cadlag lemma} and the inductive assumption we have:
\begin{align*}
&\mathbf{E}\{\sup_{s\in(\rho,\tau]}I^2_{\alpha}[g]_{\rho}^{s}(x_1,...,x_{l-k(\alpha)})\}=
\mathbf{E}\Big\{\sup_{s\in(\rho,\tau]}\int_{\rho}^{s}\int_{B}I_{\alpha-}[g]_{\rho}^{u-}(x_1,...,x_{l-k(\alpha)+1})\tilde{N}(du,dx_{l-k(\alpha)+1})\Big\}^2\\[2ex]
&\leq4\delta\int_{B}\mathbf{E}\Big\{\sup_{u\in(\rho,\tau]}I^2_{\alpha-}[g]_{\rho}^{u-}(x_1,...,x_{l-k(\alpha)+1})\Big\}\nu(dx_{l-k(\alpha)+1})\\[2ex]
&\leq4\delta \
\delta^{l(\alpha-)+s(\alpha-)-1}4^{w(\alpha-)+\tilde{n}(\alpha-)}\Big\{2(4+\delta\nu(B^{'}))\Big\}^{n(\alpha-)}\cdot\\[2ex]
&\cdot\int\limits_{B}\int_{\rho}^{\tau}\mathbf{E}\int\limits_{B_1^{\alpha-}}\int\limits_{B_2^{\alpha-}}...\int\limits_{B_{k(\alpha-)}^{\alpha-}}g^2(u,x_1,x_2,...x_{l})\nu(dx_{l})
\ \nu(dx_{l-1}) \ ... \ \nu(dx_{l-k(\alpha-)+1}) \ du \ \nu(dx_{l-k(\alpha)+1})\\[2ex]
&=\delta^{l(\alpha)+s(\alpha)-1}4^{w(\alpha)+\tilde{n}(\alpha)}\Big\{2(4+\delta\nu(B^{'}))\Big\}^{n(\alpha)}\cdot\\[2ex]
&\cdot\int_{\rho}^{\tau}\mathbf{E}\int\limits_{B_1^{\alpha}}\int\limits_{B_2^{\alpha}}...\int\limits_{B_{k(\alpha)}^{\alpha}}g^2(u,x_1,x_2,...x_{l})\nu(dx_{l})
\ \nu(dx_{l-1}) \ ... \ \nu(dx_{l-k(\alpha)+1}) \ du.
\end{align*}
\noindent $d) \ \alpha_{l(\alpha)}=3$; By (\ref{L17}), Lemma
\ref{cadlag lemma} and the inductive assumption we have:
\begin{align*}
&\mathbf{E}\{\sup_{s\in(\rho,\tau]}I^2_{\alpha}[g]_{\rho}^{s}(x_1,...,x_{l-k(\alpha)})\}=
\mathbf{E}\Big\{\sup_{s\in(\rho,\tau]}\int_{\rho}^{s}\int_{B^{'}}I_{\alpha-}[g]_{\rho}^{u-}(x_1,...,x_{l-k(\alpha)+1})N(du,dx_{l-k(\alpha)+1})\Big\}^2\\[2ex]
&\leq2(4\delta+\delta^2\nu(B^{'}))\int_{B^{'}}\mathbf{E}\Big\{\sup_{u\in(\rho,\tau]}I^2_{\alpha-}[g]_{\rho}^{u-}(x_1,...,x_{l-k(\alpha)+1})\Big\}\nu(dx_{l-k(\alpha)+1})\\[2ex]
&\leq2(4\delta+\delta^2\nu(B^{'})) \
\delta^{l(\alpha-)+s(\alpha-)-1}4^{w(\alpha-)+\tilde{n}(\alpha-)}\Big\{2(4+\delta\nu(B^{'}))\Big\}^{n(\alpha-)}\cdot\\[2ex]
&\cdot\int\limits_{B^{'}}\int_{\rho}^{\tau}\mathbf{E}\int\limits_{B_1^{\alpha-}}\int\limits_{B_2^{\alpha-}}...\int\limits_{B_{k(\alpha-)}^{\alpha-}}g^2(u,x_1,x_2,...x_{l})\nu(dx_{l})
\ \nu(dx_{l-1}) \ ... \ \nu(dx_{l-k(\alpha-)+1}) \ du \ \nu(dx_{l-k(\alpha)+1})\\[2ex]
&=\delta^{l(\alpha)+s(\alpha)-1}4^{w(\alpha)+\tilde{n}(\alpha)}\Big\{2(4+\delta\nu(B^{'}))\Big\}^{n(\alpha)}\cdot\\[2ex]
&\hspace{10ex}\cdot\int_{\rho}^{\tau}\mathbf{E}\int\limits_{B_1^{\alpha}}\int\limits_{B_2^{\alpha}}...\int\limits_{B_{k(\alpha)}^{\alpha}}g^2(u,x_1,x_2,...x_{l})\nu(dx_l)
\ \nu(dx_{l-1}) \ ... \ \nu(dx_{l-k(\alpha)+1}) \ du.
\end{align*}
\hfill$\square$\\

\noindent For any multiindex $\alpha\neq v$ and a process
$g=g(s,x_1,...,x_{k(\alpha)})$ we define two auxiliary functionals:
\begin{align}
F^{\alpha}_t[g]&:=\mathbf{E}\sup_{s\in[0,t]}\left(\sum_{i=0}^{n_s-1}I_{\alpha}[g]_{\tau_i}^{\tau_{i+1}}+
I_{\alpha}[g]_{\tau_{n_s}}^{s}\right)^2 \label{F alpha},\\[2ex]
G_{\rho,\tau}^{\alpha}[g]&:=\mathbf{E}\sup_{s\in[\rho,\tau]}\int\limits_{B_1^{\alpha}}\int\limits_{B_2^{\alpha}}...\int\limits_{B_{k(\alpha)}^{\alpha}}g^2(s,x_1,...,x_{k(\alpha)})\nu(dx_{k(\alpha)})...\nu(dx_1).
\label{G alpha}
\end{align}

\begin{lem}\label{lemat3}
For any multiindex $\alpha\neq\nu$ and a process $g$ s.t.
$G^{\alpha}_{0,t}[g]<\infty$ we have the following inequality:
\[F^{\alpha}_{t}[g]\leq
\begin{cases}
\ t \ \delta^{l(\alpha)+s(\alpha)-2} \ \int_{0}^{t}G^{\alpha}_{0,u}[g] \ du \ &\text{if} \quad l(\alpha)=s(\alpha) \\[2ex]
\ C(\alpha,t) \ \delta^{l(\alpha)+s(\alpha)-1} \
\int_{0}^{t}G^{\alpha}_{0,u}[g] \ du \ &\text{if} \quad
l(\alpha)\neq s(\alpha)\ .
\end{cases}
\]
\end{lem}
{\bf Proof:} We consider several cases:
\begin{enumerate}[a)]
\item
 $l(\alpha)=s(\alpha)$,
\item
\{ $w(\alpha)>0$ or  $\tilde{n}(\alpha)>0$\} and
 \begin{enumerate}[b1)]
 \item $\alpha_{l(\alpha)}=0$,
 \item $\alpha_{l(\alpha)}=1$,
 \item $\alpha_{l(\alpha)}=2$,
 \item $\alpha_{l(\alpha)}=3$,
 \end{enumerate}
\item
 $n(\alpha)>0$ and $w(\alpha)=\tilde{n}(\alpha)=0$.
\end{enumerate}
a) \ $l(\alpha)=s(\alpha)$\\
By the Schwarz inequality and Lemma \ref{lemat2} we have:
\begin{align*}
&F^\alpha_t[g]=\mathbf{E}\sup_{s\in[0,t]}\Big(\int_{0}^{s}I_{\alpha-}[g]_{\tau_{n_u}}^{u-}du\Big)^2\leq
\mathbf{E}\sup_{s\in[0,t]}\Big[s\int_{0}^{s}I^2_{\alpha-}[g]_{\tau_{n_u}}^{u-}du\Big]\leq
t\int_{0}^{t}\mathbf{E}\Big(I^2_{\alpha-}[g]_{\tau_{n_u}}^{u-}\Big)du\\[2ex]
&\leq
t\int_{0}^{t}\mathbf{E}\sup_{s\in(\tau_{n_u},u]}\Big(I^2_{\alpha-}[g]_{\tau_{n_u}}^{s-}\Big)du
\leq t\int_{0}^{t}\delta^{l(\alpha-)+s(\alpha-)-1}
\int_{\tau_{n_u}}^{u}\mathbf{E}g^2(s)ds \ du\\[2ex]
&\leq t\int_{0}^{t}\delta^{l(\alpha-)+s(\alpha-)-1}
\int_{\tau_{n_u}}^{u}\mathbf{E}\sup_{w\in[\tau_{n_u},s]}g^2(w)ds \
du \leq t\int_{0}^{t}\delta^{l(\alpha-)+s(\alpha-)-1}
\int_{\tau_{n_u}}^{u}G^{\alpha}_{\tau_{n_u},s}[g]ds \ du\\[2ex]
&\leq t\int_{0}^{t}\delta^{l(\alpha-)+s(\alpha-)-1} \delta
G^{\alpha}_{\tau_{n_u},u}[g] du \leq t\delta^{l(\alpha)+s(\alpha)-2}
\int_{0}^{t}G^{\alpha}_{0,u}[g] du.
\end{align*}
b1) \{ $w(\alpha)>0$ or  $\tilde{n}(\alpha)>0$\} and
$\alpha_{l(\alpha)}=0$ \\
\noindent The following inequality holds:
\begin{align*}
F^{\alpha}_{t}[g]\leq 2 \
\mathbf{E}\sup_{s\in[0,t]}\left(\sum_{i=0}^{n_s-1}I_{\alpha}[g]_{\tau_i}^{\tau_{i+1}}
\right)^2+2 \ \mathbf{E}\sup_{s\in[0,t]}\left(
I_{\alpha}[g]_{\tau_{n_s}}^{s}\right)^2.
\end{align*}
Notice that the process
$\sum_{i=0}^{n_s-1}I_{\alpha}[g]_{\tau_i}^{\tau_{i+1}}$ is a
martingale because it contains integral with respect to the Wiener
process or with respect to the compensated Poisson measure. First
let us consider the first sum.
\begin{align*}
&\mathbf{E}\sup_{s\in[0,t]}\left(\sum_{i=0}^{n_s-1}I_{\alpha}[g]_{\tau_i}^{\tau_{i+1}}
\right)^2 \leq 4
\sup_{s\in[0,t]}\mathbf{E}\left(\sum_{i=0}^{n_s-1}I_{\alpha}[g]_{\tau_i}^{\tau_{i+1}}
\right)^2 \\[2ex]
&=4
\sup_{s\in[0,t]}\mathbf{E}\left(\sum_{i=0}^{n_s-2}I_{\alpha}[g]_{\tau_i}^{\tau_{i+1}}+I_{\alpha}[g]_{\tau_{n_s-1}}^{\tau_{n_s}}
\right)^2 \\[2ex]
&=4
\sup_{s\in[0,t]}\mathbf{E}\left(\left(\sum_{i=0}^{n_s-2}I_{\alpha}[g]_{\tau_i}^{\tau_{i+1}}\right)^2+
2\sum_{i=0}^{n_s-2}I_{\alpha}[g]_{\tau_i}^{\tau_{i+1}}I_{\alpha}[g]_{\tau_{n_s-1}}^{\tau_{n_s}}+
I^2_{\alpha}[g]_{\tau_{n_s-1}}^{\tau_{n_s}} \right)\\[2ex]
 &\leq 4
\sup_{s\in[0,t]}\left\{\mathbf{E}\left(\sum_{i=0}^{n_s-2}I_{\alpha}[g]_{\tau_i}^{\tau_{i+1}}\right)^2+
\mathbf{E}\sup_{u\in(\tau_{n_s-1},\tau_{n_s}]}I^2_{\alpha}[g]_{\tau_{n_s-1}}^{u}\right\}\\[2ex]
&\leq 4 \sup_{s\in[0,t]}\Big\{\mathbf{E}\left(
\sum_{i=0}^{n_s-2}I_{\alpha}[g]_{\tau_i}^{\tau_{i+1}}\right)^2+
\delta^{l(\alpha)+s(\alpha)-1}4^{w(\alpha)+\tilde{n}(\alpha)}\Big\{2(4+\delta\nu(B^{'}))\Big\}^{n(\alpha)}\cdot\\[2ex]
&\hspace{10ex}\cdot\int_{\tau_{n_s-1}}^{\tau_{n_s}}\mathbf{E}\int\limits_{B_1^{\alpha}}\int\limits_{B_2^{\alpha}}...\int\limits_{B_{k(\alpha)}^{\alpha}}g^2(u,x_1,x_2,...x_{k(\alpha)})\nu(dx_{k(\alpha)})
\ \nu(dx_{k(\alpha)-1}) \ ... \ \nu(dx_1) du \Big\}\\[2ex]
\end{align*}
\begin{align*}
&\leq 4
\sup_{s\in[0,t]}\Bigg\{\mathbf{E}\left(\sum_{i=0}^{n_s-2}I_{\alpha}[g]_{\tau_i}^{\tau_{i+1}}\right)^2+
\delta^{l(\alpha)+s(\alpha)-1}4^{w(\alpha)+\tilde{n}(\alpha)}\Big\{2(4+\delta\nu(B^{'}))\Big\}^{n(\alpha)}\cdot\\[2ex]
&\hspace{70ex}\cdot\int_{\tau_{n_s-1}}^{\tau_{n_s}}G^{\alpha}_{\tau_{n_s-1},u}du
\Bigg\}\\[2ex]
&\leq 4
\sup_{s\in[0,t]}\Bigg\{\mathbf{E}\left(\sum_{i=0}^{n_s-3}I_{\alpha}[g]_{\tau_i}^{\tau_{i+1}}\right)^2+
\delta^{l(\alpha)+s(\alpha)-1}4^{w(\alpha)+\tilde{n}(\alpha)}\Big\{2(4+\delta\nu(B^{'}))\Big\}^{n(\alpha)}\cdot\\[2ex]
&\hspace{45ex}\cdot\left(\int_{\tau_{n_s-2}}^{\tau_{n_s-1}}G^{\alpha}_{\tau_{n_s-2},u}du+\int_{\tau_{n_s-1}}^{\tau_{n_s}}G^{\alpha}_{\tau_{n_s-1},u}du\right)
\Bigg\}\\[2ex]
&\leq 4 \sup_{s\in[0,t]}\Bigg\{
\delta^{l(\alpha)+s(\alpha)-1}4^{w(\alpha)+\tilde{n}(\alpha)}\Big\{2(4+\delta\nu(B^{'}))\Big\}^{n(\alpha)}\cdot\\[2ex]
&\hspace{30ex}\cdot\left(\int_{\tau_{0}}^{\tau_{1}}G^{\alpha}_{\tau_{0},u}du+\int_{\tau_{1}}^{\tau_{2}}G^{\alpha}_{\tau_{1},u}du+...+\int_{\tau_{n_s-1}}^{\tau_{n_s}}G^{\alpha}_{\tau_{n_s-1},u}du\right)
\Bigg\}\\[2ex]
&\leq\delta^{l(\alpha)+s(\alpha)-1}4^{w(\alpha)+\tilde{n}(\alpha)+1}\Big\{2(4+\delta\nu(B^{'}))\Big\}^{n(\alpha)}\int_{0}^{t}G^{\alpha}_{0,u}du.
\end{align*}
For the second sum we have the following inequalities:
\begin{align*}
&\mathbf{E}\sup_{s\in[0,t]}\left(
I_{\alpha}[g]_{\tau_{n_s}}^{s}\right)^2\leq\delta \
\mathbf{E}\sup_{s\in[0,t]}\int_{\tau_{n_s}}^{s}I^{2}_{\alpha-}[g]_{\tau_{n_s}}^{u-}du
\leq \delta \
\mathbf{E}\sup_{s\in[0,t]}\int_{\tau_{n_s}}^{s}\sup_{w\in(\tau_{n_s},u]}I^{2}_{\alpha-}[g]_{\tau_{n_s}}^{w-}du \\[2ex]
&\leq \delta \
\mathbf{E}\int_{0}^{t}\sup_{w\in(\tau_{n_u},u]}I^{2}_{\alpha-}[g]_{\tau_{n_u}}^{w-}du=
\delta \
\int_{0}^{t}\mathbf{E}\sup_{w\in(\tau_{n_u},u]}I^{2}_{\alpha-}[g]_{\tau_{n_u}}^{w-}du\\[2ex]
&\leq\delta \ \int_{0}^{t}
\delta^{l(\alpha-)+s(\alpha-)-1}4^{w(\alpha-)+\tilde{n}(\alpha-)}\Big\{2(4+\delta\nu(B^{'}))\Big\}^{n(\alpha-)}\cdot\\[2ex]
&\cdot\int_{\tau_{n_u}}^{u}\mathbf{E}\int\limits_{B_1^{\alpha}}\int\limits_{B_2^{\alpha}}...\int\limits_{B_{k(\alpha)}^{\alpha}}g^2(w,x_1,x_2,...x_{k(\alpha)})\nu(dx_{k(\alpha)})
\ \nu(dx_{k(\alpha)-1}) \ ... \ \nu(dx_1) dw \ du\\[2ex]
&\leq \delta \
\delta^{l(\alpha-)+s(\alpha-)-1}4^{w(\alpha-)+\tilde{n}(\alpha-)}\Big\{2(4+\delta\nu(B^{'}))\Big\}^{n(\alpha-)}
\int_{0}^{t}\int_{\tau_{n_u}}^{u}G^{\alpha}_{\tau_{n_u},w}[g] \ dw \
du\\[2ex]
&\leq \delta^2 \
\delta^{l(\alpha-)+s(\alpha-)-1}4^{w(\alpha-)+\tilde{n}(\alpha-)}\Big\{2(4+\delta\nu(B^{'}))\Big\}^{n(\alpha-)}
\int_{0}^{t}G^{\alpha}_{\tau_{n_u},u}[g] \ du\\[2ex]
&=\delta^{l(\alpha)+s(\alpha)-1}4^{w(\alpha)+\tilde{n}(\alpha)}\Big\{2(4+\delta\nu(B^{'}))\Big\}^{n(\alpha)}
\int_{0}^{t}G^{\alpha}_{0,u}[g] \ du.
\end{align*}
Finally we obtain:
\begin{gather*}
F_{t}^{\alpha}[g]\leq
2\cdot4^{w(\alpha)+\tilde{n}(\alpha)}\Big\{2(4+\delta\nu(B^{'}))\Big\}^{n(\alpha)}
\delta^{l(\alpha)+s(\alpha)-1} \int_{0}^{t}G^{\alpha}_{0,u}[g] \ du.
\end{gather*}

\noindent b2) \{ $w(\alpha)>0$ or  $\tilde{n}(\alpha)>0$\} and $\alpha_{l(\alpha)}=1$\\
\noindent By Doob's inequality, the isometric formula for Wiener
integrals and Lemma \ref{lemat2} we obtain:
\begin{align*}
&F_{t}^{\alpha}[g]=\mathbf{E}\sup_{s\leq
t}\left(\int_{0}^{s}I_{\alpha-}[g]_{\tau_{n_u}}^{u-}dW_u
\right)^2\leq 4 \sup_{s\leq
t}\mathbf{E}\left(\int_{0}^{s}I_{\alpha-}[g]_{\tau_{n_u}}^{u-}dW_u
\right)^2\\[2ex]
&=4\sup_{s\leq t}\mathbf{E}
\int_{0}^{s}I^2_{\alpha-}[g]_{\tau_{n_u}}^{u-}du= 4
\int_{0}^{t}\mathbf{E}\left(I^2_{\alpha-}[g]_{\tau_{n_u}}^{u-}\right)du\leq
4
\int_{0}^{t}\mathbf{E}\sup_{w\in(\tau_{n_u},u]}I^2_{\alpha-}[g]_{\tau_{n_u}}^{w-}du\\[2ex]
&\leq4\int_{0}^{t}\delta^{l(\alpha-)+s(\alpha-)-1}4^{w(\alpha-)+\tilde{n}(\alpha-)}\Big\{2(4+\delta\nu(B^{'}))\Big\}^{n(\alpha-)}
\int_{\tau_{n_u}}^{u}G^{\alpha-}_{\tau_{n_u},s}[g]\ ds \ du\\[2ex]
&\leq4\delta^{l(\alpha-)+s(\alpha-)-1}4^{w(\alpha-)+\tilde{n}(\alpha-)}\Big\{2(4+\delta\nu(B^{'}))\Big\}^{n(\alpha-)}
\delta\int_{0}^{t}G^{\alpha-}_{\tau_{n_u},u}[g] \ du\\[2ex]
&=\delta^{l(\alpha)+s(\alpha)-1}4^{w(\alpha)+\tilde{n}(\alpha)}\Big\{2(4+\delta\nu(B^{'}))\Big\}^{n(\alpha)}
\int_{0}^{t}G^{\alpha}_{0,u}[g] \ du.
\end{align*}

\noindent b3) \{ $w(\alpha)>0$ or  $\tilde{n}(\alpha)>0$\} and $\alpha_{l(\alpha)}=2$\\
By Doob's inequality, the isometric formula for integrals with
respect to the compensated Poisson measure and Lemma \ref{lemat2} we
obtain:
\begin{align*}
&F^{\alpha}_{t}[g]=\mathbf{E}\sup_{s\leq
t}\left(\int_{0}^{s}\int_{B}I_{\alpha-}[g]_{\tau_{n_u}}^{u-}(x_1)\tilde{N}(du,dx_1)\right)^2\leq
4\sup_{s\leq
t}\mathbf{E}\left(\int_{0}^{s}\int_{B}I_{\alpha-}[g]_{\tau_{n_u}}^{u-}(x_1)\tilde{N}(du,dx_1)\right)^2\\[2ex]
&=4\sup_{s\leq
t}\mathbf{E}\left(\int_{0}^{s}\int_{B}I^2_{\alpha-}[g]_{\tau_{n_u}}^{u-}(x_1)\nu(dx_1)du\right)
=4\mathbf{E}\left(\int_{0}^{t}\int_{B}I^2_{\alpha-}[g]_{\tau_{n_u}}^{u-}(x_1)\nu(dx_1)du\right)\\[2ex]
&\leq
4\left(\int\limits_{0}^{t}\int\limits_{B}\mathbf{E}\sup_{w\in(\tau_{n_u},u]}I^2_{\alpha-}[g]_{\tau_{n_u}}^{w-}(x_1)\nu(dx_1)du\right)\\[2ex]
&\leq 4 \int\limits_{0}^{t}\int\limits_{B}
\delta^{l(\alpha-)+s(\alpha-)-1}4^{w(\alpha-)+\tilde{n}(\alpha-)}\Big\{2(4+\delta\nu(B^{'}))\Big\}^{n(\alpha-)}\cdot\\[2ex]
&\cdot\int_{\tau_{n_u}}^{u}\mathbf{E}\int\limits_{B_1^{\alpha-}}\int\limits_{B_2^{\alpha-}}...\int\limits_{B_{k(\alpha-)}^{\alpha-}}g^2(w,x_1,x_2,...x_{k(\alpha)})\nu(dx_{k(\alpha)})
\ \nu(dx_{k(\alpha)-1}) \ ... \ \nu(dx_2) dw \ \nu(dx_1) du\\[2ex]
&\leq 4
\delta^{l(\alpha-)+s(\alpha-)-1}4^{w(\alpha-)+\tilde{n}(\alpha-)}\Big\{2(4+\delta\nu(B^{'}))\Big\}^{n(\alpha-)}\cdot\\[2ex]
&\cdot\int_{0}^{t}\mathbf{E}\int_{\tau_{n_u}}^{u}\int\limits_{B_1^{\alpha}}\int\limits_{B_2^{\alpha}}...\int\limits_{B_{k(\alpha)}^{\alpha}}g^2(w,x_1,x_2,...x_{k(\alpha)})\nu(dx_{k(\alpha)})
\ \nu(dx_{k(\alpha)-1}) \ ... \ \nu(dx_2)\nu(dx_1) dw   du\\[2ex]
\end{align*}
\begin{align*}
&\leq 4
\delta^{l(\alpha-)+s(\alpha-)-1}4^{w(\alpha-)+\tilde{n}(\alpha-)}\Big\{2(4+\delta\nu(B^{'}))\Big\}^{n(\alpha-)}\cdot\\[2ex]
&\cdot\int_{0}^{t}\mathbf{E}\int_{\tau_{n_u}}^{u}\sup_{s\in[\tau_{n_u},w]}\int\limits_{B_1^{\alpha}}\int\limits_{B_2^{\alpha}}...\int\limits_{B_{k(\alpha)}^{\alpha}}g^2(s,x_1,...x_{k(\alpha)})\nu(dx_{k(\alpha)})
\ \nu(dx_{k(\alpha)-1}) \ ... \ \nu(dx_2)\nu(dx_1) dw   du\\[2ex]
&\leq 4
\delta^{l(\alpha-)+s(\alpha-)-1}4^{w(\alpha-)+\tilde{n}(\alpha-)}\Big\{2(4+\delta\nu(B^{'}))\Big\}^{n(\alpha-)}\cdot\\[2ex]
&\cdot\int_{0}^{t}\mathbf{E}\delta\sup_{s\in[\tau_{n_u},u]}\int\limits_{B_1^{\alpha}}\int\limits_{B_2^{\alpha}}...\int\limits_{B_{k(\alpha)}^{\alpha}}g^2(s,x_1,x_2,...x_{k(\alpha)})\nu(dx_{k(\alpha)})
\ \nu(dx_{k(\alpha)-1}) \ ... \ \nu(dx_2)\nu(dx_1) du\\[2ex]
&=4\delta \
\delta^{l(\alpha-)+s(\alpha-)-1}4^{w(\alpha-)+\tilde{n}(\alpha-)}\Big\{2(4+\delta\nu(B^{'}))\Big\}^{n(\alpha-)}
\int_{0}^{t}G^{\alpha}_{\tau_{n_u},u}[g] du\\[2ex]
&\leq
\delta^{l(\alpha)+s(\alpha)-1}4^{w(\alpha)+\tilde{n}(\alpha)}\Big\{2(4+\delta\nu(B^{'}))\Big\}^{n(\alpha)}
\int_{0}^{t}G^{\alpha}_{0,u}[g] du.\\[2ex]
\end{align*}

\noindent
b4) \{ $w(\alpha)>0$ or  $\tilde{n}(\alpha)>0$\} and $\alpha_{l(\alpha)}=3$\\
We have the following inequality:
\begin{align*}
&F^{\alpha}_{t}[g]=\mathbf{E}\sup_{s\leq
t}\left(\int_{0}^{s}\int_{B^{'}}I_{\alpha-}[g]_{\tau_{n_u}}^{u-}(x_1)N(du,dx_1)\right)^2\\[2ex]
&= \mathbf{E}\sup_{s\leq
t}\left(\int_{0}^{s}\int_{B^{'}}I_{\alpha-}[g]_{\tau_{n_u}}^{u-}(x_1)\tilde{N}(du,dx_1)+
\int_{0}^{s}\int_{B^{'}}I_{\alpha-}[g]_{\tau_{n_u}}^{u-}(x_1)\nu(dx_1)du\right)^2\\[2ex]
&\leq 2 \ \mathbf{E}\sup_{s\leq
t}\left(\int_{0}^{s}\int_{B^{'}}I_{\alpha-}[g]_{\tau_{n_u}}^{u-}(x_1)\tilde{N}(du,dx_1)\right)^2+
2 \ \mathbf{E}\sup_{s\leq
t}\left(\int_{0}^{s}\int_{B^{'}}I_{\alpha-}[g]_{\tau_{n_u}}^{u-}(x_1)\nu(dx_1)du\right)^2.
\end{align*}
\noindent The first term is bounded as in the case (b3). For the
second term we have the following inequalities:
\begin{align*}
&\mathbf{E}\sup_{s\leq
t}\left(\int_{0}^{s}\int_{B^{'}}I_{\alpha-}[g]_{\tau_{n_u}}^{u-}(x_1)\nu(dx_1)du\right)^2\\[2ex]
&\leq \mathbf{E}\sup_{s\leq
t}\left(\int_{0}^{s}\int_{B^{'}}\mathbf{1} \
\nu(dx_1)du\cdot\int_{0}^{s}\int_{B^{'}}I^2_{\alpha-}[g]_{\tau_{n_u}}^{u-}(x_1)\nu(dx_1)du\right)\\[2ex]
&\leq
\delta\nu(B^{'})\int\limits_{0}^{t}\int\limits_{B^{'}}\mathbf{E}\left(I^2_{\alpha-}[g]_{\tau_{n_u}}^{u-}(x_1)\right)\nu(dx_1)du
\leq
\delta\nu(B^{'})\int\limits_{0}^{t}\int\limits_{B^{'}}\mathbf{E}\sup_{w\in(\tau_{n_u},u]}I^2_{\alpha-}[g]_{\tau_{n_u}}^{w-}(x_1)\nu(dx_1)du\\[2ex]
\end{align*}
\begin{align*}
&\leq\delta\nu(B^{'})\int_{0}^{t}\int_{B^{'}}\delta^{l(\alpha-)+s(\alpha-)-1}4^{w(\alpha-)+\tilde{n}(\alpha-)}\Big\{2(4+\delta\nu(B^{'}))\Big\}^{n(\alpha-)}\cdot\\[2ex]
&\cdot\int_{\tau_{n_u}}^{u}\mathbf{E}\int\limits_{B_1^{\alpha-}}\int\limits_{B_2^{\alpha-}}...\int\limits_{B_{k(\alpha-)}^{\alpha-}}g^2(w,x_1,x_2,...x_{k(\alpha)})\nu(dx_{k(\alpha)})
\ \nu(dx_{k(\alpha)-1}) \ ... \ \nu(dx_2) dw \ \nu(dx_1)du\\[2ex]
&=\delta\nu(B^{'})\delta^{l(\alpha-)+s(\alpha-)-1}4^{w(\alpha-)+\tilde{n}(\alpha-)}\Big\{2(4+\delta\nu(B^{'}))\Big\}^{n(\alpha-)}\cdot\\[2ex]
&\cdot\int_{0}^{t}\mathbf{E}\int_{\tau_{n_u}}^{u}\int\limits_{B_1^{\alpha}}\int\limits_{B_2^{\alpha}}...\int\limits_{B_{k(\alpha)}^{\alpha}}g^2(w,x_1,x_2,...x_{k(\alpha)})\nu(dx_{k(\alpha)})
\ \nu(dx_{k(\alpha)-1}) \ ... \ \nu(dx_1) dw du...
\end{align*}
and omitting identical operations as in (b3) we obtain:
\begin{align*}
...&\leq\delta^2\nu(B^{'})\delta^{l(\alpha-)+s(\alpha-)-1}4^{w(\alpha-)+\tilde{n}(\alpha-)}\Big\{2(4+\delta\nu(B^{'}))\Big\}^{n(\alpha-)}
\int_{0}^{t}G^{\alpha}_{0,u}[g]du\\[2ex]
&=\delta^{l(\alpha)+s(\alpha)-1}\delta\nu(B^{'})4^{w(\alpha)+\tilde{n}(\alpha)}\Big\{2(4+\delta\nu(B^{'}))\Big\}^{n(\alpha)-1}\int_{0}^{t}G^{\alpha}_{0,u}[g]du.
\end{align*}
Finally, for this case we have:
\begin{align*}
F^{\alpha}_{t}[g]&\leq 2 \
\delta^{l(\alpha)+s(\alpha)-1}4^{w(\alpha)+\tilde{n}(\alpha)}\Big\{2(4+\delta\nu(B^{'}))\Big\}^{n(\alpha)}
\int_{0}^{t}G^{\alpha}_{0,u}[g]
du\\[2ex]
&+2 \ \delta^{l(\alpha)+s(\alpha)-1}\delta\nu(B^{'})4^{w(\alpha)+\tilde{n}(\alpha)}\Big\{2(4+\delta\nu(B^{'}))\Big\}^{n(\alpha)-1}\int_{0}^{t}G^{\alpha}_{0,u}[g]du\\[2ex]
&\leq 2 \
\delta^{l(\alpha)+s(\alpha)-1}4^{w(\alpha)+\tilde{n}(\alpha)}\Big\{2(4+\delta\nu(B^{'}))\Big\}^{n(\alpha)-1}\Big\{8+3\delta\nu(B^{'})\Big\}\int_{0}^{t}G^{\alpha}_{0,u}[g]du.
\end{align*}

\noindent c) $n(\alpha)>0$ and $w(\alpha)=\tilde{n}(\alpha)=0$\\
\noindent In this case the multiindex $\alpha$ consists of $0$ and
$3$ only. If $\alpha_{l(\alpha)}=3$ then the desired inequality
follows from (b4). In opposite case let us denote
$r(\alpha):=\max\{i:\alpha_{i}=3\}$. For simplicity of exposition we
show the case when $r(\alpha)=l(\alpha)-1$. The idea for other cases
is exactly the same. We have the following inequality:
\begin{align*}
F_{t}^{\alpha}[g]\leq & 2 \ \mathbf{E}\sup_{s\leq
t}\bigg(\int_{0}^{s}
\ \int\limits_{\tau_{n_u}}^{u-}\int\limits_{B^{'}}I_{\alpha--}[g]_{\tau_{n_u}}^{w-}(x_1)\tilde{N}(dx_1,dw) \ du\bigg)^2\\[2ex]
+&2 \ \mathbf{E}\sup_{s\leq t}\bigg(\int_{0}^{s} \
\int\limits_{\tau_{n_u}}^{u-}\int\limits_{B^{'}}I_{\alpha--}[g]_{\tau_{n_u}}^{w-}(x_1)\nu(dx_1)dw
\ du\bigg)^2.
\end{align*}
Calculations for the first term in the sum above are covered by
(b1). Applying the Schwarz inequality and Lemma \ref{lemat2} for the
second term we obtain:
\begin{align*}
&\mathbf{E}\sup_{s\leq t}\bigg(\int_{0}^{s} \
\int\limits_{\tau_{n_u}}^{u-}\int\limits_{B^{'}}I_{\alpha--}[g]_{\tau_{n_u}}^{w-}(x_1)\nu(dx_1)dw
\ du\bigg)^2\leq t
\mathbf{E}\int_{0}^{t}\bigg(\int\limits_{\tau_{n_u}}^{u-}\int\limits_{B^{'}}I_{\alpha--}[g]_{\tau_{n_u}}^{w-}(x_1)\nu(dx_1)dw\bigg)^2du\\[2ex]
\end{align*}
\begin{align*}
&\leq t
\int_{0}^{t}\mathbf{E}\sup_{s\in(\tau_{n_u},u]}\bigg(\int\limits_{\tau_{n_u}}^{s-}\int\limits_{B^{'}}I_{\alpha--}[g]_{\tau_{n_u}}^{w-}(x_1)\nu(dx_1)dw\bigg)^2du\\[2ex]
&\leq t \delta \nu(B^{'})\int_{0}^{t} \
\int\limits_{\tau_{n_u}}^{u}\int\limits_{B^{'}}\mathbf{E}I_{\alpha--}^{2}[g]_{\tau_{n_u}}^{w-}(x_1)\nu(dx_1)dw
\ du\\[2ex]
&\leq t \delta \nu(B^{'})\int_{0}^{t} \
\int\limits_{\tau_{n_u}}^{u}\int\limits_{B^{'}}\mathbf{E}\sup_{s\in(\tau_{n_u},w]}I_{\alpha--}^{2}[g]_{\tau_{n_u}}^{s-}(x_1)\nu(dx_1)dw
\ du\\[2ex]
&\leq t \delta
\nu(B^{'})\delta^{l(\alpha--)+s(\alpha--)-1}\Big\{2(4+\delta\nu(B^{'}))\Big\}^{n(\alpha--)}\cdot\\[2ex]
&\cdot\int_{0}^{t} \ \int\limits_{\tau_{n_u}}^{u}\int\limits_{B^{'}}
\
\int\limits_{\tau_{n_u}}^{w}\mathbf{E}\int\limits_{B_1^{\alpha--}}...\int\limits_{B_{k(\alpha--)}^{\alpha--}}
g^2(s,x_1,...,x_{k(\alpha)})\nu(dx_{k(\alpha)})...\nu(dx_2)ds \
\nu(dx_1) \ dw\ du
\end{align*}
\begin{align*}
&\leq t
\delta^{l(\alpha--)+s(\alpha--)}\nu(B^{'})\Big\{2(4+\delta\nu(B^{'}))\Big\}^{n(\alpha--)}\cdot\\[2ex]
&\hspace{10ex}\int_{0}^{t} \ \int\limits_{\tau_{n_u}}^{u} \
\int\limits_{\tau_{n_u}}^{w}\mathbf{E}\int\limits_{B_1^{\alpha}}...\int\limits_{B_{k(\alpha)}^{\alpha}}
g^2(s,x_1,...,x_{k(\alpha)})\nu(dx_{k(\alpha)})...\nu(dx_1) \ ds \ dw\ du\\[2ex]
&\leq t
\delta^{l(\alpha--)+s(\alpha--)}\nu(B^{'})\Big\{2(4+\delta\nu(B^{'}))\Big\}^{n(\alpha--)}\int_{0}^{t}
\ \int\limits_{\tau_{n_u}}^{u} \
\delta G^{\alpha}_{\tau_{n_u},w}\ dw\ du\\[2ex]
&\leq t
\delta^{l(\alpha--)+s(\alpha--)}\nu(B^{'})\Big\{2(4+\delta\nu(B^{'}))\Big\}^{n(\alpha--)}\int_{0}^{t}
\delta^2 G^{\alpha}_{\tau_{n_u},u} du\\[2ex]
&\leq t
\delta^{l(\alpha)+s(\alpha)-1}\nu(B^{'})\Big\{2(4+\delta\nu(B^{'}))\Big\}^{n(\alpha)-1}\int_{0}^{t}
G^{\alpha}_{0,u} du.
\end{align*}
Finally we have:
\begin{align*}
F_{t}^{\alpha}[g]\leq\delta^{l(\alpha)+s(\alpha)-1}\left\{2(4+\delta\nu(B^{'}))\right\}^{n(\alpha)-1}\left\{8+2t\nu(B^{'})\right\}\int_{0}^{t}G^{\alpha}_{0,u}[g]du.
\end{align*}
\hfill$\square$

\begin{prop}\label{prop 2}
Let $\mathcal{A}$ be any hierarchical set. If for each
$\alpha\in\mathcal{A}$ the condition:
\begin{gather}\label{zalozenie w prop 2}
\int_{B_1^{\alpha}}\int_{B_2^{\alpha}}...\int_{B_{k(\alpha)}^{\alpha}}\mid
f_{\alpha}(y,x_1,x_2,...,x_{k(\alpha)})\mid^2\nu(dx_{k(\alpha)})...\nu(dx_1)\leq
L_{\alpha} ( 1+ y^2 )
\end{gather}
holds, then the approximation $Y^{\delta}$ given by
(\ref{wzor_postac_aprox.}) satisfies:
\begin{gather*}
\mathbf{E}\sup_{0\leq s\leq T}|Y^{\delta}_{s}|^2\leq
C_3(1+|Y_{0}^{\delta}|^2) \quad \forall \ t\in[0,T],
\end{gather*}
where $C_3\geq 0$.
\end{prop}
{\bf Proof:} Due to (\ref{wzor_postac_aprox.}) we write the
approximation in the following form
\begin{align*}
Y^\delta_s=Y^\delta_0&+\sum_{\alpha\in\mathcal{A}\setminus\{v\}}\left(\sum_{i=0}^{n_s-1}
 I_{\alpha}[f_{\alpha}(Y^\delta_{\tau_i},x_1,...,x_{k(\alpha)})]_{\tau_i}^{\tau_{i+1}}+ I_{\alpha}[f_{\alpha}(Y^\delta_{\tau_{n_s}},x_1,...,x_{k(\alpha)})]_{\tau_{n_s}}^{s}
\right).
\end{align*}
By Lemma \ref{lemat3} and assumption (\ref{zalozenie w prop 2}) we
have the following inequalities:
\begin{align*}
&\mathbf{E}\sup_{s\leq
t}|Y^{\delta}_{s}|^2\leq\sharp(\mathcal{A})\bigg\{|Y_{0}^{\delta}|^2+
\sum_{\alpha\in\mathcal{A}} \mathbf{E}\sup_{0\leq s\leq
T}\Big(\sum_{i=0}^{n_s-1}
 I_{\alpha}[f_{\alpha}(Y^\delta_{\tau_i},x_1,...,x_{k(\alpha)})]_{\tau_i}^{\tau_{i+1}}\\[2ex]
&\hspace{50ex}+I_{\alpha}[f_{\alpha}(Y^\delta_{\tau_{n_s}},x_1,...,x_{k(\alpha)})]_{\tau_{n_s}}^{s}
\Big)^2 \bigg\}\\[2ex]
&\leq\sharp(\mathcal{A})\bigg\{|Y_{0}^{\delta}|^2+\sum_{\alpha\in\mathcal{A}}\int_{0}^{t}\mathbf{E}\sup_{u\in[0,s]}
\int\limits_{B_1^{\alpha}}\int\limits_{B_2^{\alpha}}...\int\limits_{B_{k(\alpha)}^{\alpha}}
f_{\alpha}^2(Y^{\delta}_{\tau_{n_u}},x_1,...,x_{k(\alpha)})\nu(dx_{k(\alpha)})...\nu(dx_1)
\ ds \bigg\}
\end{align*}
\begin{align*}
&\leq\sharp(\mathcal{A})\bigg\{|Y_{0}^{\delta}|^2+\sum_{\alpha\in\mathcal{A}}\int_{0}^{t}\mathbf{E}\sup_{u\in[0,s]}
\int\limits_{B_1^{\alpha}}\int\limits_{B_2^{\alpha}}...\int\limits_{B_{k(\alpha)}^{\alpha}}
f_{\alpha}^2(Y^{\delta}_{u},x_1,...,x_{k(\alpha)})\nu(dx_{k(\alpha)})...\nu(dx_1)
\ ds \bigg\}\\[2ex]
&\leq\sharp(\mathcal{A})\bigg\{|Y_{0}^{\delta}|^2+\sum_{\alpha\in\mathcal{A}}\int_{0}^{t}\mathbf{E}\sup_{u\in[0,s]}
L_{\alpha}(1+|Y_{u}^{\delta}|^2)
\ ds \bigg\}\\[2ex]
&\leq\sharp(\mathcal{A})\bigg\{|Y_{0}^{\delta}|^2+T\sum_{\alpha\in\mathcal{A}}L_{\alpha}
+\sum_{\alpha\in\mathcal{A}}L_{\alpha}\int_{0}^{t}\mathbf{E}\sup_{u\leq
s}|Y_{u}^{\delta}|^2ds\bigg\}.
\end{align*}
By applying the Gronwall lemma \ref{Gronwal lemma} we obtain the
required result. \hfill$\square$\\

\noindent Now we are ready to present the main result's proof.\\

\noindent {\bf Proof of Theorem \ref{glowne twierdzenie}:} We write
the solution $Y$ of (\ref{rownanie_glowne}) and its approximation
$Y^\delta$ in the forms:
\begin{align}\label{rozwiazanie lamane}
Y_s=Y_0&+\sum_{\alpha\in\mathcal{A}_\gamma\setminus\{v\}}\left(\sum_{i=0}^{n_s-1}
 I_{\alpha}[f_{\alpha}(Y_{\tau_i},x_1,...,x_{k(\alpha)})]_{\tau_i}^{\tau_{i+1}}+ I_{\alpha}[f_{\alpha}(Y_{\tau_{n_s}},x_1,...,x_{k(\alpha)})]_{\tau_{n_s}}^{s}
\right)\\[2ex]\nonumber
&+\sum_{\alpha\in\mathcal{B}(\mathcal{A}_\gamma)}\left(\sum_{i=0}^{n_s-1}
 I_{\alpha}[f_{\alpha}(Y_{\bullet-},x_1,...,x_{k(\alpha)})]_{\tau_i}^{\tau_{i+1}}+ I_{\alpha}[f_{\alpha}(Y_{\bullet-},x_1,...,x_{k(\alpha)})]_{\tau_{n_s}}^{s}
\right),
\end{align}
\begin{align*}
Y^\delta_s=Y^\delta_0&+\sum_{\alpha\in\mathcal{A}_\gamma\setminus\{v\}}\left(\sum_{i=0}^{n_s-1}
 I_{\alpha}[f_{\alpha}(Y^\delta_{\tau_i},x_1,...,x_{k(\alpha)})]_{\tau_i}^{\tau_{i+1}}+ I_{\alpha}[f_{\alpha}(Y^\delta_{\tau_{n_s}},x_1,...,x_{k(\alpha)})]_{\tau_{n_s}}^{s}
\right).
\end{align*}
Due to Proposition \ref{sup drugi moment rozw } and Proposition
\ref{prop 2} the error of the approximation\\
$Z_t:=\mathbf{E}\sup_{s\leq t}\mid Y_s-Y^\delta_s\mid^2$ is finite
and satisfies the  inequality:
\begin{align}\label{oszacowanie bledu}
Z_t\leq D_1(\gamma)\left(\mid
Y_0-Y^\delta_0\mid^2+\sum_{\alpha\in\mathcal{A}_\gamma\setminus\{v\}}R^\alpha_t+\sum_{\alpha\in\mathcal{B}(\mathcal{A}_\gamma)}U^{\alpha}_{t}\right),
\end{align}
where
\begin{align}
D_1(\gamma)&={\sharp\{\mathcal{A}_\gamma\cup\mathcal{B}(\mathcal{A}_\gamma)\}},\nonumber\\[2ex]
\label{R alpha} R^{\alpha}_{t}&=\mathbf{E}\sup_{s\leq
t}\bigg(\sum_{i=0}^{n_s-1}
 I_{\alpha}\Big[f_{\alpha}(Y_{\tau_i},x_1,...,x_{k(\alpha)})-f_{\alpha}(Y^\delta_{\tau_i},x_1,...,x_{k(\alpha)})\Big]_{\tau_i}^{\tau_{i+1}}\\[2ex]
 &\hspace{20ex}+I_{\alpha}\Big[f_{\alpha}(Y_{\tau_{n_s}},x_1,...,x_{k(\alpha)})-f_{\alpha}(Y^\delta_{\tau_{n_s}},x_1,...,x_{k(\alpha)})\Big]_{\tau_{n_s}}^{s}
 \bigg)^2,\nonumber\\[2ex]\label{U alpha}
U^{\alpha}_{t}&=\mathbf{E}\sup_{s\leq t}\bigg(\sum_{i=0}^{n_s-1}
 I_{\alpha}[f_{\alpha}(Y_{\bullet-},x_1,...,x_{k(\alpha)})]_{\tau_i}^{\tau_{i+1}}+ I_{\alpha}[f_{\alpha}(Y_{\bullet-},x_1,...,x_{k(\alpha)})]_{\tau_{n_s}}^{s}
\bigg)^2.
\end{align}
Let us denote $D(\alpha,T):=\sup_{t\in[0,T]}\max\{t,C(\alpha,t)\}$
where $C(\alpha,t)$ is a constant from Lemma (\ref{lemat3}). Since
$\delta^{l(\alpha)+s(\alpha)-1}<\delta^{l(\alpha)+s(\alpha)-2}<1$,
by Lemma \ref{lemat3} and assumption {\bf(A3)} we have the following
inequality for any $\alpha\in\mathcal{A}_\gamma\backslash \{v\}$:
\begin{align*}
&R_{t}^{\alpha}\leq
D(\alpha,T)\int_{0}^{t}\mathbf{E}\sup_{s\leq u}\Big[f_{\alpha}(Y_{\tau_{n_s}},x_1,...,x_{k(\alpha)})-f_{\alpha}(Y^\delta_{\tau_{n_s}},x_1,...,x_{k(\alpha)})\Big]du\\[2ex]
&\leq D(\alpha,T)\int_{0}^{t}\mathbf{E}\sup_{s\leq
u}\int\limits_{B^\alpha_1}\int\limits_{B^\alpha_2}...\int\limits_{B^\alpha_{k(\alpha)}}\Big[f_{\alpha}(Y_s,x_1,...,x_{k(\alpha)})
\\[2ex]
&\hspace{40ex}-f_{\alpha}(Y^{\delta}_{s},x_1,...,x_{k(\alpha)})\Big]^2\nu(dx_{k(\alpha)})...\nu(dx_1) \ du\\[2ex]
&\leq D(\alpha,T)K_\alpha\int_{0}^{t}\mathbf{E}\sup_{s\leq u}\mid
Y_s-Y_{s}^{\delta}\mid^2 \ du = D(\alpha,T)K_\alpha\int_{0}^{t}Z_u \
du.
\end{align*}
\noindent For any $\alpha\in\mathcal{B}(\mathcal{A}_\gamma)$
inequality: $l(\alpha)+s(\alpha)-1>l(\alpha)+s(\alpha)-2\geq
2\gamma$ is satisfied. Due to this fact, assumption {\bf(A4)},
Proposition \ref{sup drugi moment rozw } and Lemma \ref{cadlag
lemma} we have the following inequalities:
\begin{align*}
U_{t}^{\alpha}&\leq
D(\alpha,T)\delta^{2\gamma}\int_{0}^{t}G^\alpha_{0,u}\Big[f_{\alpha}(Y_{\bullet-},x_1,...,x_{k(\alpha)})\Big]du\\[2ex]
&\leq D(\alpha,T)\delta^{2\gamma}\int_{0}^{t}\mathbf{E}\sup_{s\leq
u}\int\limits_{B^\alpha_1}\int\limits_{B^\alpha_2}...\int\limits_{B^\alpha_{k(\alpha)}}\Big[f_{\alpha}(Y_{s-},x_1,...,x_{k(\alpha)})\Big]^2\nu(dx_{k(\alpha)})...\nu(dx_1) \ du\\[2ex]
&\leq
D(\alpha,T)\delta^{2\gamma}L_\alpha\int_{0}^{t}\mathbf{E}\sup_{s\leq
u}\mid 1+Y_{s_{-}}^2\mid \ du \leq
D(\alpha,t)\delta^{2\gamma}L_\alpha\int_{0}^{t}( 1+C_2(1+Y_0^2)) \
du\\[2ex]
&\leq\delta^{2\gamma}D(\alpha,T)L_\alpha T ( 1+C_2(1+Y_0^2)).
\end{align*}
Finally, denoting shorter relevant constants we have:
\begin{gather*}
R_{t}^{\alpha}\leq D_2(\alpha,T)\int_{0}^{t}Z_u \ du,\qquad
U_{t}^{\alpha}\leq D_3(\alpha,T,Y_0)\delta^{2\gamma}.
\end{gather*}
Coming back to (\ref{oszacowanie bledu}) we obtain
\begin{gather}\label{oszacowanie bledu do Gronwalla}
Z_t\leq D_1(\gamma)\mid
Y_0-Y_0^{\delta}\mid^2+\tilde{D}_2(\gamma,T)\int_{0}^{t}Z_u
du+\tilde{D}_3(\gamma,T,Y_0)\delta^{2\gamma},
\end{gather}
where
$\tilde{D}_2(\gamma,T):=D_1(\gamma)\sum_{\alpha\in\mathcal{A}_\gamma\setminus
\{v\}}D_2(\alpha,T)$ and\\
$\tilde{D}_3(\gamma,T,Y_0):=D_1(\gamma)\sum_{\alpha\in\mathcal{B}(\mathcal{A}_\gamma)}D_3(\alpha,T,Y_0)$.
Applying the Gronwall lemma \ref{Gronwal lemma} to (\ref{oszacowanie
bledu do Gronwalla}) we obtain:
\begin{gather*}
Z_t\leq E_1(\gamma,T)\mid
Y_0-Y_0^\delta\mid^2+E_2(\gamma,T,Y_0)\delta^{2\gamma},
\end{gather*}
where:
\begin{align*}
E_1(\gamma,T)=D_1(\gamma)e^{\tilde{D}_2(\gamma,T)T}, \quad
E_2(\gamma,T,Y_0)=\tilde{D}_3(\gamma,Y_0,T)e^{\tilde{D}_2(\gamma,T)T}
.
\end{align*}
\hfill$\square$

\section{Infinite L\'evy measure}\label{Infinite Levy measure}
The strong approximation described by Theorem \ref{glowne
twierdzenie} can not always be easily constructed in practice even
for low order of convergence. In case when ${\nu(B)=\infty}$ the
integrals with respect to the compensated Poisson measure are
difficult to obtain even for simple integrands. In this section we
formulate alternative theorem which describes approximation with the
use of integrals which can be practically derived.\\
\noindent For a fixed $\varepsilon\in(0,1)$ we split the unit ball
$B$ into the ball $B_{\varepsilon}$ with radius $\varepsilon$ and
the disc $D_{\varepsilon}=B\backslash B_{\varepsilon}$. Our idea is
to modify the approximation given by Theorem \ref{glowne
twierdzenie} by exchanging all the integrals on unit balls with
respect to the compensated Poisson measure for integrals on discs
$D_{\varepsilon}$.\\
\noindent  For the use of this section we extend the inductive
definition of multiple stochastic integral introduced in Section
\ref{Basic definitions and facts}. To this end for any multiindex
$\alpha$ let us define a set of subscripts $\Pi(\alpha)$ consisting
of vectors
$\pi(\alpha)=(\pi_1(\alpha),\pi_2(\alpha),...,\pi_{\tilde{n}(\alpha)}(\alpha))$
of length $\tilde{n}(\alpha)$ with coordinates equal to $0$ or $1$,
i.e.
\[
\pi(\alpha)\in\Pi(\alpha) \Longleftrightarrow \
\begin{cases}
\pi_i(\alpha)=0  \ \text{or} \ \pi_i(\alpha)=1 \ \text{for} \
i=1,2,...,\tilde{n}(\alpha) \ &\text{if} \ \tilde{n}(\alpha)>0\\
v \ &\text{if} \ \tilde{n}(\alpha)=0.
\end{cases}
\]
The empty subscript $v$, i.e. the subscript of length zero is
introduced for technical reasons.\\
The subscripts for the multiindices $\alpha$ and $\alpha-$ are
related in the following way:
\[
\pi(\alpha-)= \
\begin{cases}
\pi(\alpha) \ &\text{if} \ \alpha_{l(\alpha)}=0,1,3
\\
(\pi_1(\alpha),\pi_2(\alpha),...,\pi_{\tilde{n}(\alpha)-1}(\alpha))
\ &\text{if} \ \alpha_{l(\alpha)}=2.
\end{cases}
\]

\noindent For a process $g=g(s,x_1,...,x_l)$, a multiindex $\alpha$
s.t. $k(\alpha)\leq l$ and a subscript $\pi(\alpha)\in\Pi(\alpha)$
we define the multiple integral by the induction procedure.\\
\noindent If $\tilde{n}(\alpha)=0$ then
\begin{gather*}
I^\varepsilon_{\alpha_v}[g]_{\rho}^{\tau}(x_1,...,x_l)=I_{\alpha}[g]_{\rho}^{\tau}(x_1,...,x_l).
\end{gather*}
Assume that
$I_{{\alpha-}_{\pi(\alpha-)}}^{\varepsilon}[g]_{\rho}^{\tau}(x_1,x_2,...,x_k)$
depends on $k$ parameters, where $k\leq l$. Then:
\begin{enumerate}[1)]
\item if $\alpha_{l(\alpha)}=0$ then
\begin{gather*}
I^\varepsilon_{\alpha_{\pi(\alpha)}}[g]_{\rho}^{\tau}(x_1,...,x_{k})=\int_{\rho}^{\tau}I_{{\alpha-}_{\pi(\alpha-)}}[g]_{\rho}^{s-}(x_1,...,x_{k})ds,
\end{gather*}
\item if $\alpha_{l(\alpha)}=1$ then
\begin{gather*}
I^\varepsilon_{\alpha_{\pi(\alpha)}}[g]_{\rho}^{\tau}(x_1,...,x_{k})=\int_{\rho}^{\tau}I_{\alpha-_{{\pi(\alpha-)}}}[g]_{\rho}^{s-}(x_1,...,x_{k})dW_s,
\end{gather*}
\item if $\alpha_{l(\alpha)}=2$ and ${\pi_{\tilde{n}(\alpha)}(\alpha)}=0$ and $k\geq 1$ then
\begin{gather*}
I^\varepsilon_{\alpha_{\pi(\alpha)}}[g]_{\rho}^{\tau}(x_1,...,x_{k-1})=\int_{\rho}^{\tau}\int_{B_\varepsilon}I_{\alpha-_{{\pi(\alpha-)}}}[g]_{\rho}^{s-}(x_1,...,x_{k})\tilde{N}(ds,dx_{k}),
\end{gather*}
\item if $\alpha_{l(\alpha)}=2$ and ${\pi_{\tilde{n}(\alpha)}(\alpha)}=1$ and $k\geq 1$ then
\begin{gather*}
I^\varepsilon_{\alpha_{\pi(\alpha)}}[g]_{\rho}^{\tau}(x_1,...,x_{k-1})=\int_{\rho}^{\tau}\int_{D_\varepsilon}I_{\alpha-_{{\pi(\alpha-)}}}[g]_{\rho}^{s-}(x_1,...,x_{k})\tilde{N}(ds,dx_{k}),
\end{gather*}
\item if $\alpha_{l(\alpha)}=3$ and $k\geq1$ then
\begin{gather*}
I^\varepsilon_{\alpha_{\pi(\alpha)}}[g]_{\rho}^{\tau}(x_1,...,x_{k-1})=\int_{\rho}^{\tau}\int_{B^{'}}I_{\alpha-_{{\pi(\alpha-)}}}[g]_{\rho}^{s-}(x_1,...,x_{l-k})N(ds,dx_{k}).
\end{gather*}
\end{enumerate}
In fact the last integral does not depend on $\varepsilon$,
nevertheless, we use this notation for technical reasons.

{\bf Example} Assume that $g$ is of the form $g(s,x_1,x_2,)$. Then:
\begin{align*}
I_{(212)_{(1,0)}}[g]_{0}^{t}&=\int\limits_{0}^{t}\int\limits_{B_\varepsilon}
\ \ \int\limits_{0}^{\ \ s_1-} \ \ \int\limits_{0}^{ \
 \ s_2-}\int\limits_{D_\varepsilon}g(s_3-,x_1,x_2) \ \tilde{N}(ds_3,dx_2) \
dW_{s_2} \ \tilde{N}(ds_1,dx_1).
\end{align*}
\noindent For any hierarchical set $\mathcal{A}$ let us denote by
${\mathcal{A}}^2$ a subset of  multiindices containing at least one
element equal to 2, i.e. $\alpha\in {\mathcal{A}}^2$ iff $\alpha\in
\mathcal{A} \ \text{and} \ \tilde{n}(\alpha)>0$.
\begin{rem}\label{rem o zamianie calek} Let $\varepsilon>0$. For any
$\alpha\in\mathcal{A}^2$ and a process
$g=g(s,x_1,x_2,...,x_{k(\alpha)})$ the following equality holds:
\begin{align*}
I_\alpha[g]_{\rho}^{\tau}=\sum_{\pi\in\Pi(\alpha)}I^\varepsilon_{\alpha_{\pi(\alpha)}}[g]_{\rho}^{\tau}.
\end{align*}
\end{rem}
\begin{rem}\label{uwaga dla lem 3}
If we replace in the formulas (\ref{F alpha}),(\ref{G alpha}) the
unit balls in integrals with respect to the compensated Poisson
measure by $\varepsilon$-balls or $\varepsilon$-discs, then Lemma
\ref{lemat3} remains true. As a consequence, for a process:
\begin{gather*}
Y_{s}^{\delta,\varepsilon}=\sum_{\alpha\in\mathcal{A}\backslash\mathcal{A}^2}
I_{\alpha}[f_{\alpha}(Y_{\tau_{n_s}}^{\delta,\varepsilon},x_1,...,x_{k(\alpha)})]_{\tau_{n_s}}^{s}
+\sum_{\alpha\in\mathcal{A}^2}
I^{\varepsilon}_{\alpha_{(1,1,...,1)}}[f_{\alpha}(Y_{\tau_{n_s}}^{\delta,\varepsilon},x_1,...,x_{k(\alpha)}^{})]_{\tau_{n_s}}^{s}
\end{gather*}
we obtain analogous estimation as in Proposition \ref{prop 2}, i.e.
\begin{gather*}
\mathbf{E}\sup_{s\leq t}|Y_{s}^{\delta,\varepsilon}|^2\leq
C_4(1+|Y_{0}^{\delta,\varepsilon}|^2) \quad \forall \ t\in[0,T],
\end{gather*}
where $C_4>0$, assuming that (\ref{zalozenie w prop 2}) is
satisfied.
\end{rem}

\begin{tw}\label{tw glowne z epsilon}
Assume that coefficients in equation (\ref{rownanie_glowne}) satisfy
conditions {\bf(A1),(A2)}. Let $\mathcal{A}_\gamma$ be a
hierarchical set given by (\ref{hierarchical set}) and assume that
{\bf (A3),(A4)} hold. Assume that for any
$\alpha\in\mathcal{A}^2_\gamma$ there exists a constant
$L_{\alpha}^{\varepsilon}$ such that for every $i$ s.t. $\alpha_i=2$
holds:
\begin{align}\label{war dodatkowy}
\int\limits_{B^{\alpha}_1}\int\limits_{B^{\alpha}_2}...\int\limits_{B_\varepsilon}...\int\limits_{B^{\alpha}_{k(\alpha)}}\mid
f_{\alpha}(y,x_1,x_2,...,x_{k(\alpha)})\mid^2\nu(dx_{k(\alpha)})...\nu(dx_1)&\leq
L^{\varepsilon}_{\alpha}( 1+y^2),
\end{align}
where $B_\varepsilon$ is on the position $k(\alpha)-i+1$ and
$L_{\alpha}^{\varepsilon}\underset{\varepsilon\longrightarrow
0}{\longrightarrow} 0$.\\
Then the approximation defined by the formula:
\begin{gather*}
Y_{s}^{\delta,\varepsilon}=\sum_{\alpha\in\mathcal{A}\gamma\backslash\mathcal{A}_\gamma^2}
I_{\alpha}[f_{\alpha}(Y_{\tau_{n_s}}^{\delta,\varepsilon},x_1,...,x_{k(\alpha)})]_{\tau_{n_s}}^{s}
+\sum_{\alpha\in\mathcal{A}_\gamma^2}
I^{\varepsilon}_{\alpha_{(1,1,...,1)}}[f_{\alpha}(Y_{\tau_{n_s}}^{\delta,\varepsilon},x_1,...,x_{k(\alpha)}^{})]_{\tau_{n_s}}^{s}
\end{gather*}
satisfies:
\begin{gather*}
\mathbf{E}\sup_{s\in[0,T]}\mid
Y_s-Y^{\delta,\varepsilon}_s\mid^2\leq N_1(\gamma,T)\mid
Y_0-Y^{\delta,\varepsilon}_0
\mid^2+N_2(\gamma,T,Y_0)\delta^{2\gamma}+N_3(\gamma,T,Y_0^{\delta,\varepsilon},\varepsilon),
\end{gather*}
where $N_3(\gamma,T,Y_0^{\delta,\varepsilon},\varepsilon)
\underset{\varepsilon\longrightarrow 0}{\longrightarrow} 0$.
\end{tw}
{\bf Proof:} We write the approximation in the form:
\begin{align*}
&Y^{\delta,\varepsilon}_s=\sum_{\alpha\in\mathcal{A}_\gamma}I_{\alpha}[f_{\alpha}(Y^{\delta,\varepsilon}_{\tau_{n_s}},x_1,...,x_{k(\alpha)})]_{\tau_{n_s}}^{s}-
\sum_{\alpha\in\mathcal{A}^2_\gamma}I_{\alpha}[f_{\alpha}(Y^{\delta,\varepsilon}_{\tau_{n_s}},x_1,...,x_{k(\alpha)})]_{\tau_{n_s}}^{s}\\[2ex]
&+
\sum_{\alpha\in\mathcal{A}^2_\gamma}I^\varepsilon_{\alpha_{(1,1,...,1)}}[f_{\alpha}(Y^{\delta,\varepsilon}_{\tau_{n_s}},x_1,...,x_{\tilde{n}(\alpha)+n(\alpha)})]_{\tau_{n_s}}^{s}\\[2ex]
&=Y^{\delta,\varepsilon}_0+\sum_{\alpha\in\mathcal{A}_\gamma\setminus\{v\}}\left(\sum_{i=0}^{n_s-1}
 I_{\alpha}[f_{\alpha}(Y^{\delta,\varepsilon}_{\tau_i},x_1,...,x_{k(\alpha)})]_{\tau_i}^{\tau_{i+1}}+ I_{\alpha}[f_{\alpha}(Y^{\delta,\varepsilon}_{\tau_{n_s}},x_1,...,x_{k(\alpha)})]_{\tau_{n_s}}^{s}
\right)\\[2ex]
&-\sum_{\alpha\in\mathcal{A}^2_\gamma}\left(\sum_{i=0}^{n_s-1}
 I_{\alpha}[f_{\alpha}(Y^{\delta,\varepsilon}_{\tau_i},x_1,...,x_{k(\alpha)})]_{\tau_i}^{\tau_{i+1}}+ I_{\alpha}[f_{\alpha}(Y^{\delta,\varepsilon}_{\tau_{n_s}},x_1,...,x_{k(\alpha)})]_{\tau_{n_s}}^{s}
\right)\\[2ex]
&+\sum_{\alpha\in\mathcal{A}^2_\gamma}\left(\sum_{i=0}^{n_s-1}
 I^\varepsilon_{\alpha_{(1,1,...,1)}}[f_{\alpha}(Y^{\delta,\varepsilon}_{\tau_{i}},x_1,...,x_{k(\alpha)})]_{\tau_i}^{\tau_{i+1}}+ I^\varepsilon_{\alpha_{(1,1,...,1)}}[f_{\alpha}(Y^{\delta,\varepsilon}_{\tau_{n_s}},x_1,...,x_{k(\alpha)})]_{\tau_{n_s}}^{s}
\right).
\end{align*}
By Remark \ref{uwaga dla lem 3} and Proposition \ref{sup drugi
moment rozw } the error $Z_t:=\mathbf{E}\sup_{s\leq
t}|Y_s-Y_{s}^{\delta,\varepsilon}|^2$ is finite. Taking into account
(\ref{rozwiazanie lamane}) we have:
\begin{align}\label{blad zaburzony}
Z_t\leq M_1(\gamma)\left(\mid
Y_0-Y^{\delta,\varepsilon}_0\mid^2+\sum_{\alpha\in\mathcal{A}_\gamma\setminus\{v\}}R^\alpha_t+\sum_{\alpha\in\mathcal{B}(\mathcal{A}_\gamma)}U^{\alpha}_{t}+\sum_{\alpha\in\mathcal{A}_\gamma^2}S^\alpha_t\right),
\end{align}
where
\begin{align*}
M_1(\gamma)={\sharp\{\mathcal{A}_\gamma\}+\sharp\{\mathcal{B}(\mathcal{A}_\gamma)\}+\sharp\{\mathcal{A}^2_\gamma\}},
\end{align*}
$R_{t}^{\alpha}$ is defined by (\ref{R alpha}) with $Y^\delta$
replaced by $Y^{\delta,\varepsilon}$ and $U_{t}^{\alpha}$ by (\ref{U
alpha}) and
\begin{align*}
S^{\alpha}_t=\mathbf{E}\sup_{s\leq t}\bigg(\sum_{i=0}^{n_s-1}\left(
 I_{\alpha}[f_{\alpha}(Y^{\delta,\varepsilon}_{\tau_i},x_1,...,x_{k(\alpha)})]_{\tau_i}^{\tau_{i+1}}-
 I^\varepsilon_{\alpha_{(1,1,...,1)}}[f_{\alpha}(Y^{\delta,\varepsilon}_{\tau_{i}},x_1,...,x_{k(\alpha)})]_{\tau_i}^{\tau_{i+1}}\right)\\[2ex]
+
I_{\alpha}[f_{\alpha}(Y^{\delta,\varepsilon}_{\tau_{n_s}},x_1,...,x_{k(\alpha)})]_{\tau_{n_s}}^{s}-I^\varepsilon_{\alpha_{(1,1,...,1)}}[f_{\alpha}(Y^{\delta,\varepsilon}_{\tau_{n_s}},x_1,...,x_{k(\alpha)})]_{\tau_{n_s}}^{s}\bigg)^2.
\end{align*}
Due to Remark \ref{rem o zamianie calek} we have:
\begin{align*}
S_{t}^{\alpha}&=\mathbf{E}\sup_{s\leq
t}\Bigg(\sum_{i=0}^{n_s-1}\bigg(\sum_{\pi(\alpha)\in\Pi(\alpha),
\pi(\alpha)\neq(1,1,...,1)}
 I^\varepsilon_{\alpha_{\pi(\alpha)}}[f_{\alpha}(Y^{\delta,\varepsilon}_{\tau_{i}},x_1,...,x_{k(\alpha)})]_{\tau_i}^{\tau_{i+1}}\bigg)\\[2ex]
&\hspace{20ex}+ \bigg(\sum_{\pi(\alpha)\in\Pi(\alpha),
\pi(\alpha)\neq(1,1,...,1)}
 I^\varepsilon_{\alpha_{\pi(\alpha)}}[f_{\alpha}(Y^{\delta,\varepsilon}_{\tau_{n_s}},x_1,...,x_{k(\alpha)})]_{\tau_{n_s}}^{s}\bigg)\Bigg)^2\\[2ex]
&\leq \Big({\sharp\{\Pi(\alpha)\}-1}\Big)\sum_{\pi\in\Pi(\alpha),
\pi\neq(1,1,...,1)}\mathbf{E}\sup_{s\leq t}\Bigg(\sum_{i=0}^{n_s-1}
I^\varepsilon_{\alpha_{\pi(\alpha)}}[f_{\alpha}(Y^{\delta,\varepsilon}_{\tau_{n_s}},x_1,...,x_{k(\alpha)})]_{\tau_i}^{\tau_{i+1}}\\[2ex]
&{\hspace{45ex}}+I^\varepsilon_{\alpha_{\pi(\alpha)}}[f_{\alpha}(Y^{\delta,\varepsilon}_{\tau_{n_s}},x_1,...,x_{k(\alpha)})]_{\tau_{n_s}}^{s}\Bigg)^2.
\end{align*}
\noindent In the sum above each integral contains at least one
integral on $\varepsilon$-ball. Using assumption $(\ref{war
dodatkowy})$ and Remark \ref{uwaga dla lem 3} we obtain:
\begin{align*}
S^{\alpha}_{t}&\leq
\Big({\sharp\{\Pi(\alpha)\}-1}\Big)\sum_{\pi(\alpha)\in\Pi(\alpha),
\pi\neq(1,1,...,1)}C(\alpha,t)L_{\alpha}^{\varepsilon}\int_{0}^t\mathbf{E}\sup_{s\leq
u} (1+|Y_s^{\delta,\varepsilon}|^2) du\\[2ex]
&\leq
\Big({(\sharp\{\Pi(\alpha)\}-1)}\Big)^2D(\alpha,T)L_{\alpha}^{\varepsilon}T(1+C_4(1+|Y_0^{\delta,\varepsilon}|^2))=:L^{\varepsilon}_{\alpha}\cdot
D_4(\alpha,T,Y_0^{\delta,\varepsilon}),
\end{align*}
Coming back to $(\ref{blad zaburzony})$ and using notation of
constants from the proof of Theorem $\ref{glowne twierdzenie}$ we
obtain:
\begin{gather*}
Z_t\leq M_1(\gamma)\mid
Y_0-Y_0^{\delta,\varepsilon}\mid^2+M_2(\gamma,T)\int_{0}^{t}Z_u du +
M_3(\gamma,T,Y_0)\delta^{2\gamma}+M_4(\gamma,T,Y_0^{\delta,\varepsilon},\varepsilon),
\end{gather*}
where
$M_2(\gamma,T)=M_1(\gamma)\sum_{\alpha\in\mathcal{A}_\gamma\backslash
v}D_2(\alpha,T)$,
$M_3(\gamma,T,Y_0)=M_1(\gamma)\sum_{\alpha\in\mathcal{B}(\mathcal{A}_\gamma)}D_3(\alpha,T,Y_0)$
and
$M_4(\gamma,T,Y_0^{\delta,\varepsilon},\varepsilon)=M_1(\gamma)\sum_{\alpha\in\mathcal{A}_\gamma^2}L_{\alpha}^{\varepsilon}\cdot
D_4(\alpha,t,Y_0^{\delta,\varepsilon})$. Finally, applying the
Gronwall lemma \ref{Gronwal lemma} we obtain:
\begin{gather*}
Z_t\leq N_1(\gamma,T)\mid Y_0-Y^{\delta,\varepsilon}_0
\mid^2+N_2(\gamma,T,Y_0)\delta^{2\gamma}+N_3(\gamma,T,Y_0^{\delta,\varepsilon},\varepsilon),
\end{gather*}
where $N_1(\gamma,T)=M_1(\gamma)e^{M_2(\gamma,T)T}$;
$N_2(\gamma,T,Y_0)=M_3(\gamma,T,Y_0)e^{M_2(\gamma,T)T}$;
$N_3(\gamma,T,Y_0^{\delta,\varepsilon},\varepsilon)=M_4(\gamma,T,Y_0^{\delta,\varepsilon},\varepsilon)e^{M_2(\gamma,T)T}=
e^{M_2(\gamma,T)T}
M_1(\gamma)\sum_{\alpha\in\mathcal{A}_\gamma^2}L_{\alpha}^{\varepsilon}\cdot
D_4(\alpha,t,Y_0^{\delta,\varepsilon})$.\\

\noindent \phantom{a} \hfill$\square$

\section{Examples}\label{Examples}
We present the Euler ($\gamma=\frac{1}{2}$) and Milstein
($\gamma=1$) schemes in for linear coefficients, i.e.
\begin{align*}
b(y)=by, \ \sigma(y)=\sigma y, \ F(y,x)=F y p(x), \ G(y,x)=G y q(x)
\end{align*}
where $\sigma,b,F,G$ are constants and functions $p(\cdot),q(\cdot)$
satisfy integrability conditions: $\int_{B}p^2(x)\nu(dx)<\infty$,
$\int_{B^{'}}q^2(x)\nu(dx)<\infty$. Then assumptions
{\bf(A1),\bf(A2)} are satisfied. \\ \\
For finding integrals with respect to the Poisson random measure we
use the representation of random measures, see for instance Th.
$6.5$ in \cite{PesZab}, applied to a set $E$ s.t. $\nu(E)<\infty$.
The random measure $N(\cdot,\cdot)$ can be represented as
\begin{align*}
N(t,E)=\sum_{n>0}\mathbf{1}_{[0,t]\times E}(\eta_n,\xi_n),
\end{align*}
where $\eta_n=r_1+r_2+...+r_n$ and $\{\xi_n\}, \{r_n\}$ are mutually
independent random variables  with distributions:
\begin{align*}
P(\xi_n\in A)=\frac{\nu(A\cap E)}{\nu(E)}, \ \forall
A\in\mathcal{B}(\mathbb{R}), \qquad P(r_n>s)=e^{-\nu(E)s}, \ s\geq
0.
\end{align*}
\noindent In the following constructions we assume that
$\nu(B)<\infty$ and as a consequence that
$N((\tau_i,\tau_{i+1}],B\cup B^{'})=:K(i)<\infty$. Then all the
moments of jumps generated by the Poisson random measure $N$ in the
interval $(\tau_i,\tau_{i+1}]$ form a sequence:
$\eta_1<\eta_2<...<\eta_{K(i)}$. We omit the dependence of this
sequence on $i$ to simplify notation. For the sake of clarity we use
the following notation:
\begin{align*}
\bar{\eta}_n=\min\{\eta_k: \eta_k>\eta_n \ \text{and} \ \xi_k\in
B^{'}\}\wedge\tau_{i+1},\\[2ex]
\underline{\eta}_n=\min\{\eta_k: \eta_k>\eta_n \ \text{and} \
\xi_k\in B\}\wedge\tau_{i+1}.
\end{align*}
Condition $\nu(B)<\infty$ guaranties that all the formulas below can
be practically derived. If it is not satisfied, then we apply
Theorem \ref{tw glowne z epsilon} by replacing all unit balls in the
approximation by $\varepsilon$- discs. In this case $K(i)$ and
$\underline{\eta}_n$ are defined with the use of $D_\varepsilon$
instead of $B$. Since $N((\tau_i,\tau_{i+1}],D_\varepsilon\cup
B^{'})<\infty$ the modified approximation can be calculated. We also
find the dependence of the approximation error on $\varepsilon$.\\

\noindent Notational remark: if the range of indices in the sums
below is empty, then the sum is assumed to be zero.\\

\noindent {\large \bf The Euler scheme}\\
The hierarchical set and the remainder sets are of the form
$\mathcal{A}_{\frac{1}{2}}=\{v,0,1,2,3\}$,
$\mathcal{B}(\mathcal{A}_\frac{1}{2})=\{00,10,20,30,01,11,21,31,02,12,22,32,03,13,23,33\}$.
It can be easily checked that conditions {\bf(A3), \bf(A4)} are also
satisfied. The approximation has the following form:
\begin{align*}
Y^{\delta}_{\tau_{i+1}}&=Y^{\delta}_{\tau_{i}}+I_0[f_0(Y^\delta_{\tau_i})]_{\tau_i}^{\tau_{i+1}}+
I_1[f_1(Y^\delta_{\tau_i})]_{\tau_i}^{\tau_{i+1}}+I_2[f_2(Y^\delta_{\tau_i},x)]_{\tau_i}^{\tau_{i+1}}
+I_3[f_3(Y^\delta_{\tau_i},x)]_{\tau_i}^{\tau_{i+1}}
\end{align*}
where:
\begin{align*}
I_0[f_0(Y^\delta_{\tau_i})]_{\tau_i}^{\tau_{i+1}}&=\int\limits_{\tau_{i}}^{\tau_{i+1}}b
Y^{\delta}_{\tau_{i}}ds=b
Y^{\delta}_{\tau_{i}}\bigtriangleup_i,\\[2ex]
I_1[f_1(Y^\delta_{\tau_i})]_{\tau_i}^{\tau_{i+1}}&=\int\limits_{\tau_{i}}^{\tau_{i+1}}\sigma
Y^{\delta}_{\tau_{i}}dW_s=\sigma
Y^{\delta}_{\tau_{i}}\bigtriangleup^W_i,\\[2ex]
I_2[f_2(Y^\delta_{\tau_i},x)]_{\tau_i}^{\tau_{i+1}}&=
\int\limits_{\tau_{i}}^{\tau_{i+1}}\int\limits_{B}F\cdot
Y^{\delta}_{\tau_{i}}p(x)\tilde{N}(ds,dx),\\[2ex]
\end{align*}
\begin{align*}
&=F Y^{\delta}_{\tau_{i}}\left(\sum_{n=1}^{K(i)}\mathbf{1}_{B}(\xi_n)p(\xi_n)-\bigtriangleup_i\int\limits_{B}p(x)\nu(dx)\right),\\[2ex]
I_3[f_3(Y^\delta_{\tau_i},x)]_{\tau_i}^{\tau_{i+1}}&=
\int\limits_{\tau_{i}}^{\tau_{i+1}}\int\limits_{B^{'}}G\cdot
Y^{\delta}_{\tau_{i}}q(x)N(ds,dx)=G
Y^{\delta}_{\tau_{i}}\left(\sum_{n=1}^{K(i)}\mathbf{1}_{B^{'}}(\xi_n)q(\xi_n)\right),\\[2ex]
\bigtriangleup_i&=\tau_{i+1}-\tau_{i},\qquad
\bigtriangleup^W_i=W_{\tau_{i+1}}-W_{\tau_{i}}.
\end{align*}
If $\nu(B)=\infty$ we apply Theorem $\ref{tw glowne z epsilon}$.
Notice that condition $(\ref{war dodatkowy})$ is satisfied since
\begin{gather*}
\int\limits_{B_\varepsilon}\mid
f_2(y,x)\mid^2\nu(dx)=F^2y^2\int\limits_{B_\varepsilon}p^2(x)\nu(dx)\underset{\varepsilon\rightarrow
0}{\longrightarrow} 0,
\end{gather*}
so
$L_{\alpha}^{\varepsilon}=\int\limits_{B_\varepsilon}p^2(x)\nu(dx)$.
It follows from the proof of Theorem $\ref{tw glowne z epsilon}$
that:
\begin{gather*}
N_3\Big(\frac{1}{2},T,Y_{0}^{\delta,\varepsilon},\varepsilon\Big)=
K(T,Y_{0}^{\delta,\varepsilon})\int\limits_{B_\varepsilon}p^2(x)\nu(dx),
\end{gather*}
where $K(T,Y_{0}^{\delta,\varepsilon})>0$.\\ \\

\noindent {\large \bf The Milstein scheme}\\
The hierarchical and remainder sets are of the form:
\begin{align*}
\mathcal{A}_1&=\{v,0,1,2,3,11,21,31,12,22,32,13,23,33\},\\
\mathcal{B}(\mathcal{A}_1)&=\{
00,10,20,30,
01,
02,
03,
011,111,211,311,
021,121,221,321,
031,131,231,331,\\
&\qquad 012,112,212,312,
022,122,222,322,
032,132,232,332,
013,113,213,313,\\
&\qquad 023,123,223,323,
033,133,233,333,
 \}.
\end{align*}
Assumptions ${\bf(A3),(A4)}$ are satisfied. The approximation is of
the following form:
\begin{align*}
Y^{\delta}_{\tau_{i+1}}=Y^{\delta}_{\tau_{i}}&+I_0[f_0(Y^\delta_{\tau_i})]_{\tau_i}^{\tau_{i+1}}+
I_1[f_1(Y^\delta_{\tau_i})]_{\tau_i}^{\tau_{i+1}}+I_2[f_2(Y^\delta_{\tau_i},x)]_{\tau_i}^{\tau_{i+1}}
+I_3[f_3(Y^\delta_{\tau_i},x)]_{\tau_i}^{\tau_{i+1}}\\[2ex]
&+I_{11}[f_{11}(Y^\delta_{\tau_i})]_{\tau_i}^{\tau_{i+1}}+
I_{21}[f_{21}(Y^\delta_{\tau_i})]_{\tau_i}^{\tau_{i+1}}+I_{31}[f_{31}(Y^\delta_{\tau_i},x)]_{\tau_i}^{\tau_{i+1}}\\[2ex]
&+I_{12}[f_{12}(Y^\delta_{\tau_i})]_{\tau_i}^{\tau_{i+1}}+
I_{22}[f_{22}(Y^\delta_{\tau_i})]_{\tau_i}^{\tau_{i+1}}+I_{32}[f_{32}(Y^\delta_{\tau_i},x)]_{\tau_i}^{\tau_{i+1}}\\[2ex]
&+I_{13}[f_{13}(Y^\delta_{\tau_i})]_{\tau_i}^{\tau_{i+1}}+
I_{23}[f_{23}(Y^\delta_{\tau_i})]_{\tau_i}^{\tau_{i+1}}+I_{33}[f_{33}(Y^\delta_{\tau_i},x)]_{\tau_i}^{\tau_{i+1}},
\end{align*}
where $I_0[f_0(Y^\delta_{\tau_i})]_{\tau_i}^{\tau_{i+1}},
I_1[f_1(Y^\delta_{\tau_i})]_{\tau_i}^{\tau_{i+1}},I_2[f_2(Y^\delta_{\tau_i},x)]_{\tau_i}^{\tau_{i+1}},
I_3[f_3(Y^\delta_{\tau_i},x)]_{\tau_i}^{\tau_{i+1}}$ are like in the
Euler scheme and
\begin{align*}
I_{11}[f_{11}(Y^\delta_{\tau_i})]_{\tau_i}^{\tau_{i+1}}&=\frac{1}{2}\sigma^2Y^\delta_{\tau_i}\left((\bigtriangleup^W_i)^2-\bigtriangleup_i\right),\\[2ex]
I_{13}[f_{13}(Y^\delta_{\tau_i},x)]_{\tau_i}^{\tau_{i+1}}&=G\sigma
Y_{\tau_i}\int_{\tau_i}^{\tau_{i+1}}\int_{B^{'}}q(x)(W_s-W_{\tau_i})N(ds,dx)\\[2ex]
&=G\sigma
Y_{\tau_i}\sum_{n=1}^{K(i)}q(\xi_n)(W_{\eta_n}-W_{\tau_i})\mathbf{1}_{B^{'}}(\xi_n),\\[2ex]
\end{align*}
\begin{align*}
I_{12}[f_{12}(Y^\delta_{\tau_i},x)]_{\tau_i}^{\tau_{i+1}}&=F\sigma
Y_{\tau_i}\bigg(\int_{\tau_i}^{\tau_{i+1}}\int_{B}p(x)(W_s-W_{\tau_i})N(ds,dx)\\[2ex]
&-\int_{\tau_i}^{\tau_{i+1}}\int_{B}p(x)(W_s-W_{\tau_i})\nu(dx)ds\bigg)\\[2ex]
&=F\sigma
Y_{\tau_i}\left(\sum_{n=1}^{K(i)}p(\xi_n)(W_{\eta_n}-W_{\tau_i})\mathbf{1}_{B}(\xi_n)-\int_{B}p(x)\nu(dx)\cdot\triangle^{Z}_{i}
\right),
\end{align*}
where
$\triangle_{i}^{Z}=\int_{\tau_i}^{\tau_{i+1}}(W_s-W_{\tau_i})ds$ is
a random variable with distribution
$N(0,\frac{1}{3}\triangle_{i}^{3})$, correlated with
$\triangle_{i}^{W}$, i.e.
$\mathbf{E}(\triangle_{i}^{W}\triangle_{i}^{Z})=\frac{1}{2}\triangle_{i}^{2}$.
The pair $\triangle_{i}^{W},\triangle_{i}^{Z}$ can be generated by
transformation of two independent random variables $U_1$,$U_2$ with
distributions $N(0,1)$ in the following way:
$\triangle_{i}^{W}=U_1\sqrt{\triangle_i}$,
$\triangle_{i}^{Z}=\frac{1}{2}\triangle_{i}^{\frac{3}{2}}(U_1+\frac{1}{\sqrt{3}}U_2)$,
for more details see \cite{KloPLa}.
\begin{align*}
I_{31}[f_{31}(Y^\delta_{\tau_i},x)]_{\tau_i}^{\tau_{i+1}}&=G\sigma
Y_{\tau_i}\int_{\tau_i}^{\tau_{i+1}}\int_{\tau_i}^{s-}\int_{B^{'}}q(x)N(du,dx)
dW_s\\[2ex]
&=G\sigma Y_{\tau_i}\sum_{n=1}^{K(i)}\Big\{ \sum_{0<k\leq
n}q(\xi_k)\mathbf{1}_{B^{'}}(\xi_k)\Big\}(W_{\bar{\eta}_n}-W_{\eta_n})
\end{align*}
\begin{align*}
I_{21}[f_{21}(Y^\delta_{\tau_i},x)]_{\tau_i}^{\tau_{i+1}}&=F\sigma
Y_{\tau_i}\bigg(\int_{\tau_i}^{\tau_{i+1}}\int_{\tau_i}^{s-}\int_{B}p(x)N(du,dx)
dW_s\\[2ex]
&-\int_{\tau_i}^{\tau_{i+1}}\int_{\tau_i}^{s-}\int_{B}p(x)\nu(dx)du \ dW_s\bigg)\\[2ex]
&=F\sigma Y_{\tau_i}\bigg(\sum_{n=1}^{K(i)}\Big\{ \sum_{0<k\leq
n}p(\xi_k)\mathbf{1}_{B}(\xi_k)\Big\}(W_{\underline{\eta}_n}-W_{\eta_n})\\[2ex]
&-\int_{B}p(x)\nu(dx)\cdot(\triangle^{W}_{i}\triangle_{i}-\triangle_{i}^{Z})\bigg)
\end{align*}
\begin{align*}
I_{33}[f_{33}(Y_{\tau_i},x_1,x_2)]_{\tau_i}^{\tau_{i+1}}&=G^2Y_{\tau_i}\int\limits_{\tau_i}^{\tau_{i+1}}\int\limits_{B^{'}}
\ \int\limits_{\tau_i}^{s-}\int\limits_{B^{'}}q(x_1)q(x_2)N(du,dx_2)
\ N(ds,dx_1)\\[2ex]
&=G^2 Y_{\tau_i}\sum_{n=1}^{K(i)}\Big\{\sum_{0<k<n}
q(\xi_k)\mathbf{1}_{B^{'}}(\xi_k)
\Big\}q(\xi_n)\mathbf{1}_{B^{'}}(\xi_n)
\end{align*}
\begin{align*}
I_{32}[f_{32}(Y_{\tau_i},x_1,x_2)]_{\tau_i}^{\tau_{i+1}}&=FGY_{\tau_i}\Big(\int\limits_{\tau_i}^{\tau_{i+1}}\int\limits_{B}
\
\int\limits_{\tau_i}^{s-}\int\limits_{B^{'}}p(x_1)q(x_2)N(du,dx_2)N(ds,dx_1)\\[2ex]
&-\int\limits_{\tau_i}^{\tau_{i+1}}\int\limits_{B} \
\int\limits_{\tau_i}^{s-}\int\limits_{B^{'}}p(x_1)q(x_2)N(du,dx_2)\nu(dx_1)ds\Big)\\[2ex]
&=FGY_{\tau_i}\Big(\sum_{n=1}^{K(i)}\Big\{\sum_{0<k<n}
q(\xi_k)\mathbf{1}_{B^{'}}(\xi_k) \Big\}p(\xi_n)\mathbf{1}_{B}(\xi_n)\\[2ex]
&-\int\limits_{B}p(x_1)\nu(dx_1)\cdot\sum_{n=1}^{K(i)}\Big\{
\sum_{0<k\leq
n}q(\xi_k)\mathbf{1}_{B}(\xi_k)\Big\}(\underline{\eta}_n-\eta_n)\Big)
\end{align*}
\begin{align*}
I_{23}[f_{23}(Y_{\tau_i},x_1,x_2)]_{\tau_i}^{\tau_{i+1}}&=FGY_{\tau_i}\bigg(\int\limits_{\tau_i}^{\tau_{i+1}}\int\limits_{B^{'}}
\
\int\limits_{\tau_i}^{s-}\int\limits_{B}q(x_1)p(x_2)N(du,dx_2)N(ds,dx_1)\\[2ex]
&-\int\limits_{\tau_i}^{\tau_{i+1}}\int\limits_{B^{'}} \
\int\limits_{\tau_i}^{s-}\int\limits_{B}q(x_1)p(x_2)du \ \nu(dx_2) \
N(ds,dx_1)\bigg)\\[2ex]
&=FGY_{\tau_i}\Big(\sum_{n=1}^{K(i)}\Big\{\sum_{0<k<n}
q(\xi_k)\mathbf{1}_{B}(\xi_k)
\Big\}p(\xi_n)\mathbf{1}_{B^{'}}(\xi_n)\\[2ex]
&-\int_{B}p(x_2)\nu(dx_2)\cdot\sum_{n=1}^{K(i)}(\eta_n-\tau_i)q(\xi_n)\mathbf{1}_{B^{'}}(\xi_n)\Big)
\end{align*}
\begin{align*}
I_{22}[f_{22}(Y_{\tau_i},x_1,x_2)]_{\tau_i}^{\tau_{i+1}}&=F^2Y_{\tau_i}\bigg(
\int\limits_{\tau_i}^{\tau_{i+1}}\int\limits_{B} \
\int\limits_{\tau_i}^{s-}\int\limits_{B}p(x_1)p(x_2)N(du,dx_2)N(ds,dx_1)\\[2ex]
&-\int\limits_{\tau_i}^{\tau_{i+1}}\int\limits_{B} \
\int\limits_{\tau_i}^{s-}\int\limits_{B}p(x_1)p(x_2)\nu(dx_2)du
 \ N(ds,dx_1) \\[2ex]
&-\int\limits_{\tau_i}^{\tau_{i+1}}\int\limits_{B} \
\int\limits_{\tau_i}^{s-}\int\limits_{B}p(x_1)p(x_2)
 \ N(du,dx_2) \ \nu(dx_1)ds\\[2ex]
&+\int\limits_{\tau_i}^{\tau_{i+1}}\int\limits_{B}
\int\limits_{\tau_i}^{s-}\int\limits_{B}p(x_1)p(x_2)
 \ \nu(dx_2)du \ \nu(dx_1)ds\bigg)\\[2ex]
\end{align*}
\begin{align*}
&=F^2Y_{\tau_i}\bigg(\sum_{n=1}^{K(i)}\Big\{\sum_{0<k<n}
p(\xi_k)\mathbf{1}_{B}(\xi_k)
\Big\}p(\xi_n)\mathbf{1}_{B}(\xi_n)\\[2ex]&
-\int_{B}p(x_2)\nu(dx_2)\cdot\sum_{n=1}^{K(i)}(\eta_n-\tau_i)p(\xi_n)\mathbf{1}_{B}(\xi_n)\\[2ex]
&-\int_{B}p(x_1)\nu(dx_1)\sum_{n=1}^{K(i)}\Big\{\sum_{0<k\leq
n}p(\xi_k)\mathbf{1}_{B}(\xi_k)\Big\}(\underline{\eta}_n-\eta_n)\\[2ex]
&+\frac{1}{2}\int_{B}p(x_1)\nu(dx_1)\int_{B}p(x_2)\nu(dx_2)\triangle_i^2
\bigg)
\end{align*}
It is easy to check that all of the integrals below:
\begin{gather*}
\int_{B_\varepsilon}|f_2(y,x)|^2\nu(dx),
\int_{B_\varepsilon}|f_{12}(y,x)|^2\nu(dx),
\int_{B_\varepsilon}|f_{21}(y,x)|^2\nu(dx)\\[2ex]
\int_{B_\varepsilon}\int_{B}|f_{22}(y,x_1,x_2)|^2\nu(dx_2)\nu(dx_1),
\int_{B}\int_{B_\varepsilon}|f_{22}(y,x_1,x_2)|^2\nu(dx_2)\nu(dx_1)\\[2ex]
\int_{B}\int_{B_\varepsilon}|f_{23}(y,x_1,x_2)|^2\nu(dx_2)\nu(dx_1),
\int_{B_\varepsilon}\int_{B}|f_{32}(y,x_1,x_2)|^2\nu(dx_2)\nu(dx_1),
\end{gather*}
are bounded above by $Ky^2\int_{B_\varepsilon}p^2(x)\nu(dx)$ where
$K$ is some constant, so we can assume that
$L_{\alpha}^{\varepsilon}=L^{\varepsilon}=\int_{B_{\varepsilon}}p^2(x)\nu(dx)$
for all $\alpha\in\mathcal{A}_{1}^{2}$. The part of the error of the
modified approximation connected with the procedure of
$\varepsilon$-balls cutting satisfies:
\begin{gather*}
N_3(1,T,Y_{0}^{\delta,\varepsilon},\varepsilon)\leq
K(T,Y_{0}^{\delta,\varepsilon})\int\limits_{B_\varepsilon}p^2(x)\nu(dx),
\end{gather*}
where $K(T,Y_{0}^{\delta,\varepsilon})>0$.\\

\end{document}